\documentclass[12pt,leqno]{article}
\usepackage{amssymb}
\newcommand{\dd}{{\rm \kern 3pt I\kern-9pt d}}

\newcommand{\Abar}{{\backslash\kern-8pt A}}

\title{\large WHEN AND HOW AN ERROR YIELDS A DIRICHLET FORM}
\author{\sc Nicolas Bouleau\footnote{ENPC 28 rue des Saints P\`eres, 75007 Paris; 
e-mail : {\tt bouleau@enpc.fr}}}
\date{\it Ecole des Ponts, ParisTech}
\begin{document}
\maketitle
{\small
\noindent{\bf Abstract.} We consider a random variable $Y$  
 and approximations $Y_n$, $n\in \mathbb{N}$, defined on  the same probability space  with values in the same measurable space as $Y$.
We are interested in situations where the approximations $Y_n$ allow to define a Dirichlet form in the space $L^2(\mathbb{P}_Y)$
where $\mathbb{P}_Y$ is the law of $Y$. Our approach consists in studying both biases and variances.
The article attempts to propose a general theoretical framework. It is illustrated by several examples.\\



\noindent {\it Contents

I. Bias operators

II. Examples

\hspace{.5cm} 1. Error in the Glivenko-Cantelli theorem

\hspace{.5cm} 2. Typical formulae of finite dimensional error calculus

\hspace{.5cm} 3. Conditionally Gaussian case

\hspace{.5cm} 4. Natural inaccuracy of the brownian motion

\hspace{.5cm} 5. Empirical laws and Brownian bridge

\hspace{.5cm} 6. Erroneous empirical laws and generalized Mehler type structures

\hspace{.5cm} 7. Erroneous random walk and Donsker theorem

\hspace{.5cm} 8. Approximation of the Brownian motion as centered orthogonal measure

\hspace{.5cm} 9. Approximation of a Poisson point process

\hspace{.5cm} 10. Stochastic integral

\hspace{.5cm} 11. Stochastic differential equation and Euler scheme

III. Conclusive comments.}}\\

\newpage

\noindent{\Large\sf I. Bias operators.}\\

We consider a random variable $Y$ defined on a probability space $(\Omega, {\cal A}, \mathbb{P})$ with values in a measurable
space $(E,{\cal F})$ and approximations $Y_n$, $n\in \mathbb{N}$, also defined on  $(\Omega, {\cal A}, \mathbb{P})$ with values in $(E,{\cal F})$.
In the whole study we suppose that there exist an algebra ${\cal D}$ of bounded functions
 from $E$ into $\mathbb{R}$ or $\mathbb{C}$ dense in $L^2(E,{\cal F},\mathbb{P}_Y)$ containing the constants and a sequence $(\alpha_n)_{n\in\mathbb{N}}$
of positive numbers, with which the following hypotheses are considered :
$$
(\mbox{H}1)\qquad\left\{\begin{array}{l}
\forall \varphi\in{\cal D}, \mbox{ there exists } \overline{A}[\varphi]\in L^2(E,{\cal F},\mathbb{P}_Y)\quad s.t. \quad\forall \chi\in{\cal D}\\
\lim_{n\rightarrow\infty} \alpha_n\mathbb{E}[(\varphi(Y_n)-\varphi(Y))\chi(Y)]=\mathbb{E}_Y[\overline{A}[\varphi]\chi].
\end{array}\right.
$$
the expectation $\mathbb{E}_Y$ being relative to the law $\mathbb{P}_Y$.
$$
(\mbox{H}2)\qquad\left\{\begin{array}{l}
\forall \varphi\in{\cal D}, \mbox{ there exists } \underline{A}[\varphi]\in L^2(E,{\cal F},\mathbb{P}_Y)\quad s.t. \quad\forall \chi\in{\cal D}\\
\lim_{n\rightarrow\infty} \alpha_n\mathbb{E}[(\varphi(Y)-\varphi(Y_n))\chi(Y_n)]=\mathbb{E}_Y[\underline{A}[\varphi]\chi].
\end{array}\right.
$$
$$
(\mbox{H}3)\quad\left\{\begin{array}{l}
\forall \varphi\in{\cal D}, \mbox{ there exists } \widetilde{A}[\varphi]\in L^2(E,{\cal F},\mathbb{P}_Y)\quad s.t. \quad\forall \chi\in{\cal D}\\
\lim_{n\rightarrow\infty} \alpha_n\mathbb{E}[(\varphi(Y_n)-\varphi(Y))(\chi(Y_n)-\chi(Y))]=-2\mathbb{E}_Y[\widetilde{A}[\varphi]\chi].
\end{array}\right.
$$
We first note that as soon as two of hypotheses (H1) (H2) (H3) are fulfilled (with
 the same algebra ${\cal D}$ and the same sequence $\alpha_n$), the third one follows thanks to the relation
$$\widetilde{A}=\frac{\overline{A}+\underline{A}}{2}.$$
When defined, the operator $\overline{A}$ which considers the asymptotic error from the point of view of the limit model, will be called {\it the 
theoretical bias operator}.

The operator $\underline{A}$ which considers the asymptotic error from the point of view of the approximating model will be called {\it the 
practical bias operator}.

Because of the property
$$<\widetilde{A}[\varphi],\chi>_{L^2(\mathbb{P}_Y)}=<\varphi,\widetilde{A}[\chi]>_{L^2(\mathbb{P}_Y)}$$
the operator $\widetilde{A}$ will be called {\it the symmetric bias operator}.\\

\noindent{\bf Remark 1.} Under (H1) the limit of $\alpha_n\mathbb{E}[(\varphi(Y)-\varphi(Y_n))\chi(Y_n)]$ exists and is equal to
$$\mathbb{E}_Y[\overline{A}[\chi]\varphi]-\mathbb{E}_Y[\overline{A}[\varphi\chi]].$$
The operator $\overline{A}$ with dense domain possesses an adjoint $\overline{A}^{\ast}$. If ${\cal D}\subset {\cal D}(\overline{A}^{\ast})$ then 
(H2) is satisfied and
\begin{equation}
\underline{A}[\varphi]=\overline{A}^{\ast}[\varphi]-\varphi\overline{A}^{\ast}[1]. \qquad\forall \varphi\in{\cal D}
\end{equation}
Reciprocally, if (H1) and (H2) are supposed and if $1\in{\cal D}(\overline{A}^\ast)$, the map 
$\chi\mapsto\mathbb{E}_Y[\underline{A}[\varphi]\chi]=\mathbb{E}_Y[\overline{A}[\chi]\varphi]-\mathbb{E}_Y[\overline{A}[\varphi\chi]]$ is continuous
and so is $\chi\mapsto\mathbb{E}[\overline{A}[\chi]\varphi]$ which shows ${\cal D}\subset{\cal D}\overline{A}^\ast$ and relation (1) holds. We see that
the hypothesis $1\in{\cal D}(\overline{A}^\ast)$ is rather strong, it will be not fulfilled in general. \\

\noindent{\bf Example I.1.} Let us take for $(E,{\cal F})$ a metrisable compact set with its Borel $\sigma$-field and let $(X_t)$ 
be a Feller process
 with values in $E$ and transition semi-group $(P_t)$ (A Feller process on an l.c.d. space reduces to this situation by the Alexandrov compactification
see [12] chap. XIII \S 20-21).

Let $(A_P,D_P)$ be the generator of the strongly continuous  contraction semi-group $(P_t)$ on ${\cal C}(E)$. Suppose $(P_t)$ be in duality with a strongly continuous
semi-group $(Q_t)$ with respect to a probability measure $\nu$ and let $(A_Q,D_Q)$ be the generator of $(Q_t)$ on ${\cal C}(E)$. 
Then, if there
is an algebra of bounded functions containing constants ${\cal D}\subset D_P\cap D_Q$ dense in $L^2(\nu)$, the 
approximation $X_t$ of $X_0$ satisfies hypotheses (H1) to (H3) and we have on ${\cal D}$:
$$\begin{array}{rl}
\overline{A}[\varphi]=&A_P[\varphi]\\
\underline{A}[\varphi]=&A_Q[\varphi]-\varphi A_Q[1]\\
\overline{A}^\ast[\varphi]=&\underline{A}[\varphi]+\varphi A_Q[1].
\end{array}
$$
Indeed, if $\varphi, \chi\in{\cal D}$
$$\frac{1}{t}\mathbb{E}_{\nu}[(\varphi(X_t)-\varphi(X_0))\chi(X_0)]=
\frac{1}{t}<P_t\varphi-\varphi,\chi>_\nu\quad\rightarrow\quad <A_P[\varphi],\chi>_\nu$$
$$\frac{1}{t}\mathbb{E}_{\nu}[(\varphi(X_0)-\varphi(X_t))\chi(X_t)]=
\frac{1}{t}[<Q_t\varphi-\varphi+\varphi(1-Q_t1),\chi>_\nu]\rightarrow
 <A_Q[\varphi]-\varphi A_Q[1],\chi>_\nu.$$
Hence (H1) to (H3) are fulfilled and, by theorem 1 below, the limit $\frac{1}{t}\mathbb{E}_\nu[(\varphi(X_t)-\varphi(X_0)^2]$ extends to 
a symmetric Dirichlet form on $L^2(\nu)$.\hfill$\diamond$\\

The basis of our study is the following theorem :\\

\noindent{\bf Theorem 1.}{\it Under hypothesis} (H3) {\it

a) the limit
$$\widetilde{\cal E}[\varphi,\chi]=\lim_n  \frac{\alpha_n}{2}\mathbb{E}[(\varphi(Y_n)-\varphi(Y))(\chi(Y_n)-\chi(Y)]\qquad \varphi, \chi\in{\cal D}$$
defines a closable positive bilinear form whose smallest closed extension is denoted $({\cal E},\mathbb{D})$.

b) $({\cal E},\mathbb{D})$ is a Dirichlet form

c) $({\cal E},\mathbb{D})$ admits a square field operator $\Gamma$ satisfying $\forall \varphi,\chi\in{\cal D}$
$$
\Gamma[\varphi]=\widetilde{A}[\varphi^2]-2\varphi\widetilde{A}[\varphi]$$
$$\mathbb{E}_Y[\Gamma[\varphi]\chi]=\lim_n\alpha_n\mathbb{E}[(\varphi(Y_n)-\varphi(Y))^2(\chi(Y_n)+\chi(Y))/2]$$
\indent d) $({\cal E},\mathbb{D})$ is local if and only if $\forall \varphi\in{\cal D}$
$$\lim_n \alpha_n\mathbb{E}[(\varphi(Y_n)-\varphi(Y))^4]=0.$$}

\noindent{\bf Demonstration.} a) That $(\widetilde{\cal E}, {\cal D})$ be closable comes from the Friedrich construction of the minimal selfadjoint
extension of a symmetric operator. Let us recall the argument.

By $\widetilde{\cal E}[\varphi,\chi]=-<\widetilde{A}[\varphi],\chi>_{L^2(\mathbb{P}_Y)}$ $\forall \varphi, \chi\in{\cal D}$ the form 
$(\widetilde{\cal E}, {\cal D})$ satisfies 
$$u_n\in{\cal D}, \|u_n\|\rightarrow 0\quad\Rightarrow \quad{\cal E}[u_n,v]\rightarrow 0\quad\forall v\in{\cal D}$$
and this property suffices to imply closability (cf. [13] ex. 1.1.2, [9] Chap. I ex. 1.3.4 or [5] lemma III.24).

b) In order to prove that the form $({\cal E},\mathbb{D})$ is Dirichlet we will use the following elementary property :

{\it If $K$ is a compact subset of $\mathbb{R}$, $\forall \varepsilon>0$ there exists a polynomial $p(x)$ such that

(i) $0\leq p(y)-p(x)\leq y-x\quad\forall x<y\in K$

(ii) $|p(x)-x|\leq \varepsilon\quad \forall x\in[0,1]\cap K$

(iii) $p(x)\geq -\varepsilon\quad\forall x\in K.$}

\noindent Let $(R_\lambda)_{\lambda>0}$ be the strongly continuous contraction resolvent associated with $({\cal E},\mathbb{D})$, we have to prove that 
the operators $\lambda R_\lambda$ are sub-Markov (cf. [13], [9], [23]). For that, since 
here $\lambda R_\lambda1=1,\;\forall\lambda>0$, because ${\cal E}[1,u]=0\;\forall u\in\mathbb{D}$, it is enough to show that 
$0\leq u\leq 1\;\Rightarrow\;R_\lambda u\geq 0$ and this for $\lambda\geq1$ since 
$R_\alpha=\frac{1}{\beta}\sum_{k=1}^\infty (\beta R_{\alpha+\beta})^k\quad\forall \alpha,\beta>0$
by the resolvent equation.

Let $u$ be a measurable function from $E$ into $\mathbb{R}$ s.t. $0\leq u\leq 1$, denoting as usual 
${\cal E}_\lambda$ the form ${\cal E}+\lambda\|.\|^2_{L^2}$, the functional
$$\Phi(v)={\cal E}[v]+\lambda\|v-\frac{u}{\lambda}\|^2\quad v\in\mathbb{D}$$
satisfies
$$\Phi(R_\lambda u)+{\cal E}_\lambda[R_\lambda u-v]=\Phi(v)\quad v\in\mathbb{D}.$$
It follows that $R_\lambda u$ is the unique element of $\mathbb{D}$ minimizing $\Phi$ on $\mathbb{D}$, and the sequence 
$v_n$ converges in $\mathbb{D}$ (equipped
with the norm $\sqrt{{\cal E}_1}$) to $R_\lambda u$ if and only if $\lim_n \Phi(v_n)=\Phi(R_\lambda u)$.

Since ${\cal D}$ is dense in $\mathbb{D}$, let $v_n\in{\cal D}$ be such a sequence converging to $R_\lambda u$. Let $p_n$ be a 
polynomial fulfilling
property {\it(i)} to {\it(iii)} for $\varepsilon=\frac{1}{n}$ and $K$ containing the range of the bounded function $v_n$. Let us 
put $w_n=p_n\circ v_n$. Since ${\cal D}$ is an algebra, $w_n\in{\cal D}$ and by {\it(i)} we 
have $\widetilde{{\cal E}}[w_n]\leq\widetilde{\cal E}[v_n]$. Now, by {\it(i)} and {\it (ii)} 
$|p_n\circ v_n-y|\leq |v_n-y|+\frac{1}{n}\quad\forall y\in[0,1]$
so $|w_n-\frac{u}{\lambda}|\leq|v_n-\frac{u}{\lambda}|+\frac{1}{n}$ what gives
$\|w_n-\frac{u}{\lambda}\|^2_{L^2}\leq\|v_n-\frac{u}{\lambda}\|^2_{L^2}+\frac{2}{n}\|v_n-\frac{u}{\lambda}\|_{L^2}+\frac{1}{n^2}$ and 
$$\Phi(w_n)\leq\Phi(v_n)+\frac{2}{n}\|v_n-\frac{u}{\lambda}\|+\frac{1}{n^2}.$$
The sequence $v_n$ being bounded in $L^2$, it follows that $\Phi(w_n)\rightarrow\Phi(R_\lambda u)$ and 
$w_n\rightarrow R_\lambda u$ in $\mathbb{D}$.
Taking eventually a subsequence converging a.s.  and using $w_n\geq -\frac{1}{n}$ by {\it (iii)}, we obtain $R_\lambda u\geq 0$ 
what proves the property.

c) Let us denote $(A,{\cal D}A)$ the selfadjoint operator associated with $({\cal E},\mathbb{D})$ (Friedrich's extension of 
$(\widetilde{A},{\cal D}))$,
since the algebra ${\cal D}\subset{\cal D}A$ is dense in $\mathbb{D}$, the theorem 4.2.2 of [9] applies and 
the form $({\cal E},\mathbb{D})$ admits a square field operator satisfying $\Gamma[\varphi]=
\widetilde{A}[\varphi^2]-2\varphi\widetilde{A}[\varphi]\quad\forall\varphi\in{\cal D}$ and $\Gamma[\psi]=A[\psi^2]-2\psi A[\psi]$ if $\psi\in{\cal D}A$ and 
$\psi^2\in{\cal D}A$. The formula of the statement comes from 
$$\mathbb{E}_Y[\Gamma[\varphi]\chi]=\mathbb{E}_Y[\widetilde{A}[\varphi^2]\chi]-2\mathbb{E}_Y[\widetilde{A}[\varphi]\varphi\chi]$$
coming back to the definition of $\widetilde{A}$.

d) Let us remark that for $f\in{\cal D}$ we have
$$4\widetilde{\cal E}[f^3,f]-3\widetilde{\cal E}[f^2,f^2]=\lim_n\alpha_n\mathbb{E}[(f(Y_n)-f(Y))^4].$$
So, if the assumption of the statement holds, using the fact that one dimensional contractions are continuous on $\mathbb{D}$ (cf. [])
$$4{\cal E}[u^3,u]-3{\cal E}[u^2,u^2]=0\quad\forall u\in\mathbb{D}\cap L^\infty$$ this is enough to imply that 
${\cal E}$ is local (cf [9] Chap. I \S 5). Reciprocally, if ${\cal E}$ is local, since ${\cal E}$ admits a square field operator, the functional calculus applies
(cf. [9] Chap. I \S 6) and
$$4{\cal E}[u^3,u]-3{\cal E}[u^2,u^2]=2\mathbb{E}[3u^2\Gamma[u]]-\frac{3}{2}\mathbb{E}[4u^2\Gamma[u]]=0$$
$\forall u\in\mathbb{D}\cap L^\infty$ and the condition of the statement is fulfilled.\hfill$\diamond$\\

\noindent{\bf Comment. } Considering for $Y$ a Brownian motion $B$ indexed by $[0,1]$ as random variable with values in ${\cal C}([0,1])$ and taking
 for $Y_\varepsilon$ the approximation $Y_\varepsilon=B+\sqrt{\varepsilon}W$ where $W$ is an independent standard Bronian motion, we may
apply the theorem with ${\cal D}$ the linear combinations of functions $\varphi(B)=e^{i\int_0^1f\,dB}$ with regular $f$ say ${\cal C}^1_b$.

We have with $\chi(B)=e^{i\int_0^1g\,dB}$
$$
\begin{array}{l}\mathbb{E}[(e^{i\int_0^1f\,dY_\varepsilon}-e^{i\int_0^1f\,dY})(e^{i\int_0^1g\,dY_\varepsilon}-e^{i\int_0^1g\,dY})]\qquad\qquad\\
\qquad\qquad\qquad=\mathbb{E}[e^{i\int (f+g)\,dY}]\mathbb{E}[(e^{i\sqrt{\varepsilon}\int fdW}-1)(e^{i\sqrt{\varepsilon}\int gdW}-1)]
\end{array}$$
so that
$$\lim_{\varepsilon\rightarrow 0}\frac{1}{\varepsilon}\mathbb{E}[(\varphi(Y_\varepsilon)-\varphi(Y)(\chi(Y_\varepsilon)-\chi(Y)]
=(-\int_0^1 fg\,dt)e^{-\frac{1}{2}\int_0^1(f+g)^2dt}$$
what may be written $-2<\widetilde{A}[\varphi],\chi>$ with
$$\widetilde{A}[\varphi]=e^{i\int f\,dB}[-\frac{i}{2}\int f\,dB-\frac{1}{2}\int f^2 dt]$$
as seen by an elementary calculation. Hypothesis (H3) is satisfied. The theorem yields the well known Ornstein-Uhlenbeck structure on the Wiener space (see e.g. [9] or [27]).

We can say that from a pedagogical point of view, in order to introduce the error calculus on the Wiener space (basic Malliavin calculus) and the same 
would be true for the Poisson space or the Monte Carlo space (cf. [5]), theorem 1 is a quite convenient tool. It is simpler 
than the theorem on infinite
products of Dirichlet structures used in [9] or [5]. It allows also to construct Dirichlet forms in a variety of situations as will 
show the examples below.

This being said, when studying mathematically approximations, the most usefull part of the theorem is often the easiest one to prove, i.e.
part a), because the limit form is often recognized as a classical form whose properties (Dirichlet character, square field operator, locality)
are known.\\

\noindent{\bf Remark 2.} Suppose that instead of (H3) we assume that the limit 
$\lim_n \alpha_n\mathbb{E}[(\varphi(Y_n)-\varphi(Y)^2]$ exists $\forall\varphi\in{\cal D}$. Then, denoting $2\widetilde{\cal E}[\varphi]$ this limit,
if the form $(\widetilde{\cal E},{\cal D})$ is closable, 

i) the hypothesis (H3) is equivalent to ${\cal D}\subset{\cal DA}$ where $({\cal A, DA})$
is the generator of the form $({\cal E},\mathbb{D})$ smallest closed extension of $(\widetilde{\cal E},{\cal D})$
(Indeed, $u\in{\cal DA}\Leftrightarrow f\mapsto{\cal E}[f,u]$ is continuous on ${\cal D}$, hence if ${\cal D}\subset{\cal DA}$ hypothesis (H3) is satisfied
with $\widetilde{A}=A$ on ${\cal D}$ and if (H3) is satified ${\cal D}\subset{\cal DA}$)

ii) parts b) and c) of the demonstration of the theorem apply and show that $({\cal E},\mathbb{D})$ is Dirichlet with square field operator.\\

\noindent{\bf Remark 3.} Under (H3) the condition d) of the theorem 
$\forall\varphi\in{\cal D}\quad\lim_n\alpha_n\mathbb{E}[(\varphi(Y_n)-\varphi(Y))^4]=0$ is equivalent to either of the conditions :

(j) $\exists\lambda>2\quad\lim_n\alpha_n\mathbb{E}[|\varphi(Y_n)-\varphi(Y)|^\lambda]=0\quad\forall\varphi\in{\cal D}.$

(jj) $\forall\lambda>2\quad\lim_n\alpha_n\mathbb{E}[|\varphi(Y_n)-\varphi(Y)|^\lambda]=0\quad\forall\varphi\in{\cal D}.$

\noindent Indeed, it suffices to proves that (j) implies (jj). Let $\lambda$ be such that (j) is true, then for $\mu>\lambda$
$$\alpha_n\mathbb{E}[|\varphi(Y_n)-\varphi(Y)|^\mu]\leq 
\alpha_n\mathbb{E}[|\varphi(Y_n)-\varphi(Y)|^\lambda](2\|\varphi\|_\infty)^{\mu-\lambda}\rightarrow 0$$
and for $\lambda>\mu=2+\varepsilon$
$$\alpha_n\mathbb{E}[|\varphi(Y_n)-\varphi(Y)|^\mu]\leq 
(\alpha_n\mathbb{E}|\varphi(Y_n)-\varphi(Y)|^2)^{1/2}(\alpha_n\mathbb{E}|\varphi(Y_n)-\varphi(Y)|^{2+2\varepsilon})^{1/2}$$
iterating the procedure until $2+2^k\varepsilon\geq \lambda$ gives the result.\hfill$\diamond$\\

\noindent{\bf Remark 4.} As an example where the limit form is non-local, let us consider the case where $(E,{\cal F})=(\Omega, {\cal A})$
and let $\theta$ be a measurable map from $\Omega$ into itself preserving the probability $\mathbb{P}$ and defining a strongly mixing
endomorphism $A\in{\cal A}\mapsto \theta^{-1}(A)$, i.e. such that
$$\lim_n\mathbb{E}[f\circ\theta^n.g]=\mathbb{E} f\mathbb{E} g\qquad\forall f,g\in L^2(\mathbb{P})$$
then taking for $Y$ the identity map and $Y_n=\theta^n$ with $\alpha_n=1$, hypotheses (H1) to (H3) hold on ${\cal D}=L^\infty(\mathbb{P})$
with $\overline{A}[\varphi]=\underline{A}[\varphi]=\widetilde{A}[\varphi]=-\varphi+\mathbb{E}[\varphi]$ and 
$\widetilde{\cal E}[\varphi]=var\varphi$.\hfill$\diamond$\\

We introduce now the fourth bias operator $\Abar$ defined under (H1) and (H2) on ${\cal D}$ as 
$$\Abar=\frac{1}{2}(\overline{A}-\underline{A}).$$
By $\mathbb{E}_Y[\Abar[\varphi]\chi]=\lim_n\mathbb{E}[(\varphi(Y_n)-\varphi(Y))(\chi(Y)+\chi(Y_n))/2]$ we see that $\Abar$ represent the asymptotic error
from the point of view of an exterior observator according the same weight to both the theoretical and the practical models and measuring
the error algebraically on the same axis. Because of the properties of $\Abar$ proved 
below, $\Abar$ will be called {\it the singular bias operator}.

We shall say that an operator $B$ from ${\cal D}$ into $L^2(\mathbb{P}_Y)$ is a {\it first order operator} if it satisfies
$$B[\varphi\chi]=B[\varphi]\chi+\varphi B[\chi]\qquad\forall\varphi,\chi\in{\cal D}$$

\noindent{\bf Proposition 1.} {\it Under} (H1) {\it to} (H3) {\it 

a) the}  theoretical variance $\lim_n\alpha_n\mathbb{E}[(\varphi(Y_n)-\varphi(Y))^2\psi(Y)]$ {\it and  the } practical variance
$\lim_n\alpha_n\mathbb{E}[(\varphi(Y_n)-\varphi(Y))^2\psi(Y_n)]$ {\it exist and we have $\forall\varphi,\chi,\psi\in{\cal D}$
$$
\begin{array}{c}
\lim_n\alpha_n\mathbb{E}[(\varphi(Y_n)-\varphi(Y))(\chi(Y_n)-\chi(Y))\psi(Y)]=\mathbb{E}_Y[-\underline{A}[\varphi\psi]\chi+\underline{A}[\psi]\varphi\chi
-\overline{A}[\varphi]\chi\psi]\\
\lim_n\alpha_n\mathbb{E}[(\varphi(Y_n)-\varphi(Y))(\chi(Y_n)-\chi(Y))\psi(Y_n)]=\mathbb{E}_Y[-\overline{A}[\varphi\psi]\chi+\overline{A}[\psi]\varphi\chi
-\underline{A}[\varphi]\chi\psi]
\end{array}
$$
\indent \it b) These two variances coincide if and only if $\;\Abar$ is a first order operator, and then are equal to 
$\mathbb{E}_Y[\Gamma[\varphi]\psi].$}\\

\noindent{\bf Proof.} The part a) comes directly from the definition of $\overline{A}$ and $\underline{A}$. The difference between the two expressions
in $\varphi,\chi,\psi$ is 
$$2\mathbb{E}_Y[(\Abar[\varphi\psi]-\Abar[\psi]\varphi-\Abar[\varphi]\psi)\chi]$$
and vanishes iff $\Abar$ is first order.\hfill$\diamond$\\

A sufficient condition for the equality of the theoretical and the practical variances is given by\\

\noindent{\bf Proposition 2.} {\it Under} (H1) {\it to} (H3) {\it If there is a real number $p\geq 1$ s.t.
$$\lim_n\alpha_n\mathbb{E}[(\varphi(Y_n)-\varphi(Y))^2|\psi(Y_n)-\psi(Y)|^p]=0\quad\forall\varphi,\psi\in{\cal D}$$
then $\Abar$ is first order.}\\

\noindent{\bf Proof.} Let be $\lambda\in[0,2[$. We have 
$|\alpha_n\mathbb{E}[(\varphi(Y_n)-\varphi(Y))^2(\psi(Y_n)-\psi(Y))]| $
$$\leq 
\alpha_n\mathbb{E}[|\varphi(Y_n)-\varphi(Y)|^\lambda|\varphi(Y_n)-\varphi(Y)|^{2-\lambda}|\psi(Y_n)-\psi(Y)|]$$
the case $p=1$ is obtained taking $\lambda=0$. If $p>1$ we go on with $\lambda>0$
$$\leq(\alpha_n\mathbb{E}[(\varphi(Y_n)-\varphi(Y))^2])^{\lambda/2}(\alpha_n\mathbb{E}[(\varphi(Y_n)-\varphi(Y))^2|\psi(Y_n)-\psi(Y)|^{2/(2-\lambda)}])^{(2-\lambda)/2}$$
the result follows taking $2/(2-\lambda)=p$.\hfill$\diamond$\\

In particular under (H1) to (H3), if the locality condition of theorem 1 is fulfilled then $\Abar$ is a first order operator.\\

\noindent{\bf Remark 5.} In the frequent cases where $\Abar$ is a first order operator, {\it $\widetilde{A}$ captures all the 
diffusive part
of $\overline{A}$ and of $\underline{A}.$}

Similarly, we can remark that for deterministic approximations, the operator $\widetilde{A}$ is often nought (see prop. 17 in fine). For example let us consider the 
ordinary differential equation
$$x_t=x_0+\int_0^tf(x_s)y_sds$$
approximated by the Euler scheme
$$x^n_t=x_0+\int_0^t f(x^n_{[ns]/n})y_sds$$
even if we suppose $x_0$ to be random, errors are of deterministic nature and as soon as $f$ is ${\cal C}^1$ with at most linear growth
and $\int_0^1y_s^2ds<+\infty$ we have for $\varphi,\chi\in{\cal C}^1_b$ (bounded with bounded derivative)
$$\begin{array}{rl}
n\mathbb{E}[(\varphi(x^n_t)-\varphi(x_t))\chi(x_t)]&\rightarrow\mathbb{E}[u_t\varphi^\prime(x_t)\chi(x_t)]\\
n\mathbb{E}[(\varphi(x_t)-\varphi(x_t^n))\chi(x_t^n)]&\rightarrow-\mathbb{E}[u_t\varphi^\prime(x_t)\chi(x_t)]\\
\mbox{and }\quad\mathbb{E}[(\varphi(x_t^n)-\varphi(x_t)^2]&\rightarrow 0
\end{array}
$$
where $u_t$ is given by $u_t=-\frac{1}{2}\int_0^t f^\prime(x_s)f(x_s)y^2_se^{\int_s^t f^\prime(x_\alpha)y_\alpha d\alpha}ds.$
(cf. [17] theorem 1.1). Thus $$\overline{A}[\varphi](x)=\mathbb{E}[u_t\varphi^\prime(x_t)|x_t=x]=-\underline{A}[\varphi](x)$$
and we have $\Abar=\overline{A}$ and $\widetilde{A}=0$.\hfill$\diamond$\\

Let us derive  some consequences of the only (H1) hypothesis. Under (H1) we may consider the symmetric {\it positive} bilinear form
$$e[\varphi,\chi]=-\mathbb{E}[\overline{A}[\varphi]\chi+\varphi\overline{A}[\chi]-\overline{A}[\varphi\chi]]$$
{\bf Proposition 3.} {\it Under } (H1) {\it the following conditions are equivalent

1) }(H2)

{\it 2)} (H3)

{\it 3) $(e,{\cal D})$ satisfies the following sufficient closability condition}
$$\chi_n\in{\cal D}\quad\chi_n\rightarrow 0 {\mbox{ in }}L^2\quad\Rightarrow
\quad e[\varphi,\chi_n]\rightarrow 0\quad \forall\varphi\in{\cal D}.$$
{\bf Proof.} Since $\lim_n \alpha_n\mathbb{E}[(\varphi(Y)-\varphi(Y_n))\chi(Y_n)]=
\mathbb{E}_Y[\overline{A}[\chi]\varphi-\overline{A}[\varphi\chi]]$, (H2) is equivalent
to 
$$\forall\varphi\in{\cal D}\quad\chi\mapsto\mathbb{E}_Y[\overline{A}[\chi]\varphi-\overline{A}[\varphi\chi]]\mbox{ is continuous on }{\cal D}
\mbox{ in }L^2$$
which is equivalent to 
$$\forall\varphi\in{\cal D}\quad\chi\mapsto
e[\varphi,\chi]\mbox{ continuous on }{\cal D}
\mbox{ in }L^2,$$
i.e. equivalent to this bilinear form be continuous at 0 which is the condition of the statement.\hfill$\diamond$\\

\noindent{\bf Proposition 4.} {\it Under} (H1){\it, the conditions of the preceding proposition are fulfilled if $\forall\varphi\in{\cal D}$}
$$\chi_p\in{\cal D},\quad\chi_p\rightarrow0\mbox{ in }L^2\quad
\Rightarrow\quad\lim_p\lim_n\alpha_n\mathbb{E}[(\chi_p(Y_n)-\chi_p(Y))\varphi(Y)]=0.$$
{\bf Proof.} The condition of the statement means $\chi\rightarrow\mathbb{E}[\overline{A}[\chi]\varphi]$ continuous, i.e. ${\cal D}\subset{\cal D}(\overline{A}^\ast)$
hence by remark 1 hypothesis (H2) holds.\hfill$\diamond$\\

\noindent{\bf Remark 6.} If $\forall\varphi\in{\cal D}$ the conditional expectation $\alpha_n\mathbb{E}[\varphi(Y_n)-\varphi(Y)|Y=y]$ converges
weakly in $L^2(\mathbb{P}_Y)$ then (H1) is satisfied, because the weak limit is necessarily an element of $L^2$.\\

\noindent{\bf Proposition 5.} {\it Under} (H1), {\it if the law of the pair $(Y_n,Y)$ is asymptotically symmetric in the following sense :
$$\lim_n\alpha_n\mathbb{E}[\varphi(Y_n)\psi(Y)-\varphi(Y)\psi(Y_n)]=0\quad\forall\varphi,\psi\in{\cal D}$$
then the conditions of proposition 3 are fulfilled, $\underline{A}=\overline{A}=\widetilde{A}$ and $\Abar=0$.}\\

\noindent{\bf Proof.} Taking $\psi=1$ gives  $\lim_n\alpha_n\mathbb{E}[\varphi(Y_n)\chi(Y_n)-\varphi(Y)\chi(Y)]=0$ hence 
$$\lim_n\alpha_n\mathbb{E}[(\varphi(Y)-\varphi(Y_n))\chi(Y_n)]=\lim_n\alpha_n\mathbb{E}[(\varphi(Y_n)-\varphi(Y))\chi(Y)]$$
and (H2) holds with $\underline{A}=\overline{A}$.\hfill$\diamond$\\

We now come back to the situation where the only assumption (H3) is supposed.\\

\noindent{\bf Theorem 2.} {\it Under} (H3). {\it If the form $({\cal E},\mathbb{D})$ (cf. theorem 1) is local, then the} principle of asymptotic error calculus
{\it is valid on 
$$\widetilde{\cal D}=\{F(f_1,\ldots,f_p)\;:\;f_i\in{\cal D},\;\;F\in{\cal C}^1(\mathbb{R}^p,\mathbb{R})\}$$
i.e.}
$
\lim_n\alpha_n\mathbb{E}[(F(f_1(Y_n),\ldots,f_p(Y_n))-F(f_1(Y),\ldots,f_p(Y))^2]$\hfill

\hfill$=\mathbb{E}_Y[\sum_{i,j=1}^p F^\prime_i(f_1,\ldots,f_p)F^\prime_j(f_1,\ldots,f_p)\Gamma[f_i,f_j]].$\\

\noindent{\bf Demonstration.} a) Let us first give the argument in the case $p=1$. When the form is local, by remark 3, 
$\lim_n\alpha_n\mathbb{E}[|f(Y_n)-f(Y)|^k]=0\;$ $\;\forall f\in{\cal D}$ for any integer $k\geq 3$. Let $F\in{\cal C}^1(\mathbb{R},\mathbb{R})$, writing the finite increments formula
$F(y)-F(x)=(y-x)\int_0^1 F^\prime(x+t(y-x))dt$ we have
$$\alpha_n\mathbb{E}[(F\circ f(Y_n)-F\circ f(Y))^2]=\alpha_n\mathbb{E}[(f(Y_n)-f(Y))^2(\int_0^1 F^\prime(f(Y)+t(f(Y_n)-f(Y)))dt)^2].$$
Let $P_k$ be a polynomial uniformly close to $F^\prime$ on the closed ball $B(0,\|f\|_\infty)$, considering
$$
\begin{array}{rl}
(I)&=\alpha_n\mathbb{E}[(f(Y_n)-f(Y))^2(\int_0^1 P_k(f(Y)+t(f(Y_n)-f(Y)))dt)^2]\\
&=\alpha_n\mathbb{E}[(f(Y_n)-f(Y))^2\left((P_k(f(Y))^2+\sum_{\ell=0,m=1}^N \lambda_{\ell m}(f(Y))^\ell(f(Y_m)-f(Y))^m\right)]\\
&=\alpha_n\mathbb{E}[(f(Y)-f(Y_n))^2\left((P_k(f(Y_n))^2+\sum_{\ell=0,m=1}^N \lambda_{\ell m}(f(Y_n))^\ell(f(Y)-f(Y_n))^m\right)]
\end{array}
$$
we observe that $(I)$ has same limit as
$$\alpha_n\mathbb{E}[(f(Y_n)-f(Y))^2\left((P_k(f(Y)))^2+(P_k(f(Y_n)))^2\right)/2]$$
which converges to $\mathbb{E}_Y[\Gamma[f]P_k^2(f)]$ by theorem 1.

Now, $|\alpha_n\mathbb{E}[(F\circ f(Y_n)-F\circ f(Y))^2] -\mathbb{E}_Y[F^{\prime 2}\!\circ\! f\;\; \Gamma[f]]|$
$$\begin{array}{c}
\leq \alpha_n\mathbb{E}\left[(f(Y_n)\!-\!f(Y))^2\left|\left(\int_0^1F^\prime(f(Y)\!+\!t(f(Y_n)\!-\!f(Y))dt\right)^2\right.\right.\qquad\qquad\qquad\qquad\qquad\qquad\\
\qquad\qquad\qquad\qquad\qquad\qquad\!-\!
\left.\left.\left(\int_0^1P_k(f(Y)\!+\!t(f(Y_n)\!-\!f(Y))dt\right)^2\right|\right]\\
+\left|\alpha_n\mathbb{E}\left[(f(Y_n)-f(Y))^2\left(\int_0^1P_k(f(Y)+t(f(Y_n)-f(Y))dt\right)^2\right]-\mathbb{E}_Y[P_k^2\!\circ\! f\, \Gamma[f]]\right|\\
\qquad\qquad\qquad+\left|\mathbb{E}_Y[P_k^2\!\circ\! f \,\Gamma[f]]-\mathbb{E}_Y[F^{\prime 2}\!\circ\! f \,\Gamma[f]]\right|.
\end{array}
$$
Since $\sup_n\alpha_n\mathbb{E}[(f(Y_n)-f(Y))^2]<+\infty$, the first and the last terms may be made small uniformly in $n$ by a suitable choice of $k$, as the second
term goes to zero when $n\uparrow+\infty$, the proof in complete in this case.

b) In the general case the finite increments formula writes
$$F(y_1,\ldots,y_p)-F(x_1,\ldots,x_p)=\sum_{i=1}^p(y_i-x_i)\int_0^1F^\prime_i(y_1,\ldots,y_{i-1},x_i+t(y_i-x_i),x_{i+1},\ldots,x_p)dt.$$
The local property of the form implies
\begin{equation}
\lim_n\alpha_n\mathbb{E}[\prod_{i=1}^k|f_i(Y_n)-f_i(Y)|]=0\quad\forall f_1,\ldots,f_k\in{\cal D}\quad\forall k\geq 3
\end{equation}
by H\"{o}lder inequality $\mathbb{E}|\prod_{i=1}^k X_i|\leq \prod_{i=1}^k(\mathbb{E}[|X_i|^k])^{1/k}$. Then the proof proceeds similarly, approximating the derivatives
$F^\prime_i$ by polynomials $P_{k,i}$ on the ball $B(0,\max_i\|f_i\|_\infty)$ of $\mathbb{R}^p$ using (2) and the property $\forall\varphi,\chi\in{\cal D}$
$$\lim_n\alpha_n\mathbb{E}[(f_i(Y_n)-f_i(Y))(f_j(Y_n)-f_j(Y))\left(\varphi(Y_n)\chi(Y)+\varphi(Y)\chi(Y_n)\right)/2]
=\mathbb{E}_Y[\Gamma[f_i,f_j]\varphi\chi]$$
which is consequence of theorem 1.\hfill$\diamond$\\

Let us end this part by a remark concerning the transportation of the four bias operators by image (cf. also [8]).\\

\noindent{\bf Remark 7.} Let $Y_n$ be an approximation of $Y$ satisfying (H1) to (H3) on the same dense algebra ${\cal D}$ of bounded functions
with constants. Let $\Phi$ be a map from $(E,{\cal F})$ to $(G,{\cal H})$ such that the algebra ${\cal D}_\Phi=\{u\;:\;u\circ\Phi\in{\cal D}\}$ be dense
in $L^2(\mathbb{P}_{\Phi\circ Y})$. Let us put $Z_n=\Phi\circ Y_n$ and $Z=\Phi\circ Y$. Then $Z_n$ and $Z$ satisfy (H1) to (H3) with ${\cal D}_\Phi$
and the same sequence $\alpha_n$ :
$$\begin{array}{c}
\overline{A}_\Phi[u](z)=\mathbb{E}_Y[\overline{A}[u\circ \Phi]\,|\,\Phi=z]\\
\underline{A}_\Phi[u](z)=\mathbb{E}_Y[\underline{A}[u\circ \Phi]\,|\,\Phi=z]
\end{array}
$$
and similar relations for $\widetilde{A}_\Phi$ and $\Abar_\Phi$. The Dirichlet form associated with $\widetilde{A}_\Phi$ is the image by 
$\Phi$ of the Dirichlet form associated with $\widetilde{A}$, and the properties of images of Dirichlet forms (cf. [9] chapter V) apply (square
field operator, locality). If $\Abar$ is first order, $\Abar_\Phi$ is first order.\\

\noindent{\Large\sf II. Examples.}\\

\noindent{\bf II.0. Preliminary example.} This is not stricly speaking an example but a part of Dirichlet forms theory itself. Let $({\cal E},\mathbb{D})$ be a 
Dirichlet form on the Hilbert space $L^2(E,{\cal F},m)$ where $m$ is a probability measure and let $(P_t)$ be the strongly continuous contraction
semi-group associated with $({\cal E},\mathbb{D})$.

Let us suppose that the quasi-regularity assumption is fulfilled so that we may construct a Markov process $Y_t$ with $P_t$ as 
transition semi-group (cf. [23] chapter IV \S3), and let us suppose also that the domain ${\cal D}A$ of the generator 
$(A,{\cal D}A)$ contains an algebra 
${\cal D}$ of bounded functions with constants dense in $L^2$. Then for $f\in{\cal D}$, the approximate forms
$${\cal E}_t[f]=\frac{1}{t}<f-P_t f,f>_{L^2(m)}=\frac{1}{2t}\mathbb{E}_m[(f(Y_0)-f(Y_t))^2]$$
do converge  (increasingly) when $t\downarrow 0$ to ${\cal E}[f]=-<Af,f>$ (cf. [13], [9] or [5]). Hence hypothesis (H3)
is fulfilled. The form $({\cal E},\mathbb{D})$ is an extension of that one provided by theorem 1.

Here, as easily seen, we have 
$$\overline{A}[f]=\underline{A}[f]=\widetilde{A}[f]=A[f]\qquad\forall f\in{\cal D}$$ and the operator $\Abar$ vanishes. The 
above properties of Dirichlet forms hold either for local or non-local forms. Since $\frac{1}{2t}\uparrow+\infty$ we see that the hypothesis (H3)
may be satisfied with $\alpha_n\uparrow+\infty$ the limit form being nevertheless non-local (cf. e) of theorem 1).\\

\noindent{\bf Note.} {\it In the whole article} the positive symmetric bilinear forms of real functions ${\cal E}[f,g]$ are extended
to complex functions, not as Hermitian forms, but as bilinear forms with the same symbol ${\cal E}$, in other words
$${\cal E}[f_1+if_2,g_1+ig_2]={\cal E}[f_1,f_2]+i{\cal E}[f_1,g_2]+i{\cal E}[f_2,g_1]-{\cal E}[f_2,g_2].$$

\noindent{\bf II.1.  Error in the Glivenko-Cantelli theorem. } Let us begin with a simple one-dimensional example related to the Glivenko-Cantelli theorem. Let $X$ be a real random 
variable with continuous distribution function $F$ anf let $X_k$ be i.i.d. copies of $X$.

If we put $$Y_n=\frac{1}{n}\sum_{k=1}^n 1_{\{X_k\leq X\}}$$
and $Y=F(X)$, we have $Y_n\rightarrow Y$ a.s. and we may study $Y_n$ as approximation of $Y$. Thinking for instance $X$ and the $X_k$'s to 
be simulated by the inversion method shows that the pair $(Y_n,Y)$ has the same law as $(\frac{1}{n}\sum_{k=1}^n1_{U_k\leq U},U)$ where 
$U$ is uniformly distributed on $[0,1]$ and the $U_k$ are i.i.d. copies of $U$. Since our framework involves only the joint law of $Y_n$ and $Y$
we may work with $(U_k,U)$ instead of $(X_k,X)$. Choosing $\alpha_n=n$ and ${\cal D}={\cal L}\{x\mapsto e^{2i\pi px},\;p\in\mathbb{Z}\}$, we obtain
that hypotheses (H1) to (H3) are fulfilled with
$$
\begin{array}{l}
\overline{A}[\varphi](y)=\frac{y-y^2}{2}\varphi^{\prime\prime}(y)\\
\underline{A}[\varphi](y)=\frac{y-y^2}{2}\varphi^{\prime\prime}(y)+(1-2y)\varphi^\prime(y)\\
\widetilde{A}[\varphi](y)=\frac{y-y^2}{2}\varphi^{\prime\prime}(y)+\frac{1-2y}{2}\varphi^\prime(y)\\
\widetilde{\cal E}[\varphi]=-<\widetilde{A}\varphi,\varphi>=\int_0^1\frac{y-y^2}{2}\varphi^{\prime 2}(y)\,dy\\
\Abar[\varphi](y)=(y-\frac{1}{2})\varphi^\prime(y)
\end{array}
$$
the theoretical and practical variances coincide and $\mathbb{D}=\{ f\in L^2[0,1]\,:\,f^\prime \mbox{ in distribution} \break\mbox{ sense }\in L^1_{loc}(]0,1[)
\mbox{ and } y\mapsto\sqrt{y-y^2}f^\prime(y)\in L^2(dy)\}$. Let us give some indication on the proof of hypothesis (H1) for instance. Denoting
$\mathbb{E}_y$ the conditional law given $Y=y$, we have to study 
$$n\int_0^1\mathbb{E}_y[(e^{2i\pi pY_n}-e^{2i\pi py})e^{2i\pi qy}]dy
=n\int_0^1((e^{2i\pi p/n}y+1-y)^n-e^{2i\pi py})e^{2i\pi qy}\,dy$$
which may be expanded as
$$=n\int_0^1\left(\exp\{-\frac{y-y^2}{2n}(2\pi p)^2\}(1+\varepsilon(n,y))-1\right)e^{2i\pi(p+q)y}\,dy$$
where $\varepsilon(n,y)$ goes to zero uniformly in $y$ when $n\rightarrow\infty$. Using $e^{-\lambda}-1=-\lambda\int_0^1e^{-t\lambda}dt$, the dominated
convergence theorem applies and the limit is 
$$\int_0^1-\frac{y-y^2}{2}(2\pi p)^2e^{2i\pi(p+q)y}\,dy=\mathbb{E}_Y[\frac{Y-Y^2}{2}\varphi^{\prime\prime}(Y)\chi(Y)]$$
for $\varphi(y)=\exp 2i\pi py$ and $\chi(y)=\exp 2i\pi qy$.\\

\noindent{\bf II.2. Typical formulae of finite dimensional error calculus.}\\

{\bf II.2.a. }Let us consider a triplet of real random variables $(Y,Z,T)$ and a real random variable $G$ independent of $(Y,Z,T)$ centered with variance one.
We are interested in the approximation $Y_\varepsilon$ of $Y$ given by
\begin{equation}
Y_\varepsilon=Y+\varepsilon Z+\sqrt{\varepsilon}TG.
\end{equation}
In the multidimensional case, $Y$ is with values in $\mathbb{R}^p$ as $Z$, $T$ is a $p\!\times\!q$-matrix and $G$ is independent of $(Y,Z,T)$ with 
values in $\mathbb{R}^q$, centered, square integrable, such that $\mathbb{E}[G_iG_j]=\delta_{ij}.$\\

\noindent{\bf Operator $\overline{A}$.}\\

\noindent{\bf Proposition  6.} {\it If $Z$ and $T$ are square integrable, if $\varphi$ is ${\cal C}^2$ bounded with bounded derivatives of first and
second orders ($\varphi\in{\cal C}^2_b$) and if $\chi$ is bounded,
$$\frac{1}{\varepsilon}\mathbb{E}[(\varphi(Y_\varepsilon)-\varphi(Y))\chi(Y)]\rightarrow\mathbb{E}_Y[\overline{A}[\varphi]\chi]$$
where $\overline{A}[\varphi](y)=\mathbb{E}[Z|Y\!=\!y]\varphi^\prime(y)+\frac{1}{2}\mathbb{E}[T^2|Y\!=\!y]\varphi^{\prime\prime}(y)$. 

In the multidimensional case 
$$\overline{A}[\varphi](y)=\mathbb{E}[Z^t|Y\!=\!y]\nabla\varphi(y)+\frac{1}{2}\sum_{ij}\mathbb{E}[(TT^t)_{ij}|Y\!=
\!y]\varphi^{\prime\prime}_{ij}(y).$$}
{\bf Proof.} Let us give the argument with the notation of the case $q=p=1$. The Taylor-Lagrange formula applied up to second order gives
$$
\begin{array}{c}\frac{1}{\varepsilon}\mathbb{E}[(\varphi(Y_\varepsilon)-\varphi(Y))\chi(Y)]=\mathbb{E}[Z\varphi^\prime(Y)\chi(Y)]\qquad
\qquad\hspace{5.5cm}\\
\qquad\qquad\quad+\frac{1}{2}
\mathbb{E}[(\varepsilon Z^2+2\sqrt{\varepsilon}ZTG+T^2G^2)\int_0^1\int_0^1\varphi^{\prime\prime}(Y+ab(\varepsilon Z+
\sqrt{\varepsilon}TG))2adadb\,\chi(Y)]
\end{array}
$$
(note that $ZTG$ and $T^2G^2\in L^1$ because of the independence) and  this converges by dominated Lebesgue theorem
to
$\mathbb{E}[Z\varphi^\prime(Y)\chi(Y)]+\frac{1}{2}\mathbb{E}[T^2\varphi^{\prime\prime}(Y)\chi(Y)].$
\hfill$\diamond$\\

\noindent{\bf Quadratic form and operator $\widetilde{A}$.}\\

\noindent{\bf Proposition 7. }{\it  If $Z$ and $T$ are square integrable, if $\varphi$ and $\chi$ are ${\cal C}^1_b$ 
$$\frac{1}{\varepsilon}\mathbb{E}[(\varphi(Y_\varepsilon)-\varphi(Y)(\chi(Y_\varepsilon)-\chi(Y)]
\rightarrow\mathbb{E}[T^2\varphi^\prime(Y)\chi^\prime(Y)]$$
and in the multidimensional case
$$\frac{1}{\varepsilon}\mathbb{E}[(\varphi(Y_\varepsilon)-\varphi(Y)(\chi(Y_\varepsilon)-\chi(Y)]
\rightarrow\mathbb{E}[(\nabla\varphi)^t(Y)TT^t\nabla\chi(Y)].$$}
{\bf Proof.} The demonstration is similar with a first order expansion.\hfill$\diamond$\\

In order to exhibit the operator $\widetilde{A}$, we must examine the conditions of an integration by parts in the preceding
limit. Let us put $\theta_{ij}(y)=\mathbb{E}[(TT^t)_{ij}|Y\!=\!y]$ so that $\mathbb{E}[(\nabla\varphi)^t(Y)TT^t\nabla\chi(Y)]=
\sum_{ij}\mathbb{E}_Y[\varphi^\prime_i\theta_{ij}\chi^\prime_i]$.\\

\noindent{\bf Proposition 8. }{\it If $Z$ and $T$ are square integrable, if for $i,j=1,\ldots,p$ the measure $\theta_{ij}\mathbb{P}_Y$ on $\mathbb{R}^p$
possesses a partial derivative in the sense of distributions $\partial_j(\theta_{ij}\mathbb{P}_Y)$ which is a bounded 
measure absolutely continuous w.r. to $\mathbb{P}_Y$, say $\rho_{ij}\mathbb{P}_Y$, then as soon as $\theta_{ij}$ and $\rho_{ij}\in L^2(\mathbb{P}_Y)$ the form $\widetilde{\cal E}[\varphi,\chi]=
\frac{1}{2}\sum_{ij}\mathbb{E}_Y[\varphi^\prime_i\theta_{ij}\chi^\prime_j]$
is closable on the algebra ${\cal D}={\cal C}^2_b$, hypotheses} (H1) {\it to} (H3) {\it are fulfilled and 
$$\widetilde{A}[\varphi]=\frac{1}{2}\sum_{ij}\theta_{ij}\varphi^{\prime\prime}_{ij}+\frac{1}{2}\sum_{ij}\rho_{ij}\varphi^\prime_j.$$}
{\bf Proof.} We have
$$\sum_{ij}\int\theta_{ij}\varphi^\prime_i\chi^\prime_j\,d\mathbb{P}_Y=
\sum_{ij}\int\theta_{ij}(\partial_j(\varphi_i^\prime\chi)-\varphi^{\prime\prime}_{ij}\chi)d\mathbb{P}_Y$$
and the equality
$$\int\theta_{ij}\partial_j(\varphi^\prime_i\chi)d\mathbb{P}_Y=-\int\varphi^\prime_i\chi\rho_{ij}d\mathbb{P}_Y$$
valid for $\varphi,\chi\in{\cal C}^\infty_K$ extends, under the assumptions of the statement, to 
$\varphi,\chi\in{\cal C}^2_b$. This yields 
$$\frac{1}{2}\sum_{ij}\mathbb{E}[\varphi^\prime_i\theta_{ij}\chi^\prime_j]=-\frac{1}{2}\int(\sum_{ij}\theta_{ij}\varphi^{\prime\prime}_{ij}
+\sum_{ij}\rho_{ij}\varphi^\prime_j)\chi\,d\mathbb{P}_Y.$$
\hfill Q.E.D.$\diamond$

The operator $\widetilde{A}$ depends only on $T$, not on $Z$. We obtain $\underline{A}$ by difference :
$$\underline{A}[\varphi]= \frac{1}{2}\sum_{ij}\theta_{ij}\varphi^{\prime\prime}_{ij}
+\sum_j(\sum_i\rho_{ij}-z_j)\varphi^\prime_j$$ where $z_j(y)=\mathbb{E}[Z_j|Y\!=\!y]$. At last, $\Abar$ is first order :
$$\Abar[\varphi]=\sum_j(z_j-\frac{1}{2}\sum_i\rho_{ij})\varphi^\prime_j.$$
{\bf Remark 9.} The results of this section II.2.a) would be identical with an approximation of the form
\begin{equation}
Y_\varepsilon=Y+\varepsilon Z+T.B_\varepsilon
\end{equation}
where $B$ is a centered Brownian motion vanishing at zero independent of $(Y,Z,T)$ 
since only the joint law of $(Y,Y_\varepsilon)$ is used.

The question is very close to the classical approach of Kolmogorov [20] to study Markov processes starting from the assumptions
$$\begin{array}{c}
\lim_{h\downarrow 0}\frac{1}{h}\mathbb{E}[(X_{t+h}-X_t)|{\cal F}_t]=b(t,X_t)\\
\\
\lim_{h\downarrow 0}\frac{1}{h}\mathbb{E}[(X_{t+h}-X_t)^2|{\cal F}_t]=a(t,X_t)
\end{array}
$$
Indeed, it is easy to see that the representation (4) occurs naturally for Ito processes and for diffusion processes given by an Ito equation :

On a filtered probability space $(\Omega,({\cal F}_t),\mathbb{P})$, let $B_t$ be an $({\cal F}_t)$-Brownian motion centered vanishing at zero and let $\xi$
be an Ito process defined by
$$\xi_t=\xi_0+\int_0^t\sigma_s\,dB_s+\int_0^tb_s\,ds$$
where the processes $\sigma$ and $b$ are adapted and continuous at zero in $L^2(\mathbb{P})$ and $L^1(\mathbb{P})$ respectively. 
Then approximating $\xi_0$ by $\xi_t$ is equivalent to approximating $\xi_0$ by
$\widetilde{\xi}_t=\xi_0+tb_0+\sigma_0B_t$ because for $\varphi\in{\cal C}^2_b$
$$\begin{array}{c}
\lim_{t\rightarrow 0}\frac{1}{t}\mathbb{E}[(\varphi(\xi_t)-\varphi(\widetilde{\xi_t})^2|{\cal F}_0]=0\\
\\
\lim_{t\rightarrow 0}\frac{1}{t}\mathbb{E}[(\varphi(\xi_t)-\varphi(\widetilde{\xi_t})|{\cal F}_0]=0
\end{array}
$$
as soon as for instance $\mathbb{E}\sigma^4_t$ is bounded in a neighborhood of zero, as seen by application of Ito formula and standard inequalities.\\

{\bf II.2.b.  Series with independent increments.}\\

Let be $$S=\sum_{n=1}^\infty \frac{X_n}{n^2}+\frac{Z_n}{n}$$ where $X_n,Z_n\in L^{2+\varepsilon}$, $Z_n$ centered, $(X_n,Z_n)$ i.i.d., we 
approximate $S$ by its partial sum $S_n=\sum_{k=1}^n\frac{X_k}{k^2}+\frac{Z_k}{k}$. 

Using Burkholder inequality, we observe that $n\mathbb{E}[|S-S_n|^{2+\varepsilon}]\rightarrow 0$ as $n\rightarrow\infty$. Thus, 
taking ${\cal D}={\cal C}^\infty_K$, we have for $\varphi,\chi\in{\cal D}$
$$\lim_n n\mathbb{E}[(\varphi(S)-\varphi(S_n))^2]=\lim_n n\mathbb{E}[(S-S_n)^2\varphi^{\prime 2}(S_n)]=\mathbb{E}[Z_1^2]\mathbb{E}[\varphi^{\prime 2}(S)]$$
$$\lim_n n\mathbb{E}[(\varphi(S)-\varphi(S_n))\chi(S_n)]=\lim_n n\mathbb{E}[(S-S_n)\varphi^\prime(S_n)\chi(S_n)+\frac{1}{2}(S-S_n)^2\varphi^{\prime\prime}(S_n)\chi(S_n)]$$
$$=\frac{1}{2}\mathbb{E}[Z_1^2]\mathbb{E}[\varphi^{\prime\prime}(S)\chi(S)]+\mathbb{E}[X_1]\mathbb{E}[\varphi^\prime(S)\chi(S)].$$
We can conclude that hypothesis (H2) is satisfied and $$\underline{A}[\varphi]=
\frac{\mathbb{E}[Z_1^2]}{2}\varphi^{\prime\prime}+ \mathbb{E}[X_1]\varphi^\prime.$$ 
Assumption (H3) holds as soon as the law of $S$ satisfies the Hamza condition ([13] p.105) and then the Dirichlet form is local.\\

{\bf II.2.c. Tails of martingales.}\\

Let us first consider the classical case of Polya's urn in its simplest configuration with two colors, one ball added each time, and an initial composition
of one white ball and one black ball.

The ratio $X_n$ of white balls after the $n$-th drawing satisfies
$$X_{n+1}(n+3)=X_n(n+1)+1_{U_{n+1}\leq X_n}$$
where $U_{n+1}$ is a random variable uniformly distributed on $[0,1]$ independent of ${\cal F}_n=\sigma(X_0,\ldots,X_n)$, i.e.
$$X_{n+1}=X_n+\frac{1}{n+3}(1_{U_{n+1}\leq X_n}-X_n).$$
Let $X$ be the (a.s. and $L^p,\;1\leq p<+\infty$) limit of the bounded martingale $X_n$, we study the approximation of $X$ by $X_n$.

We note that $\lim_n n\mathbb{E}[|X-X_n|^3]=0$ as easily seen using Burkholder inequality. Then, taking for ${\cal D}$ the 
functions of class ${\cal C}^3$ on $[0,1]$ vanishing at 0 and 1, we have $\forall \varphi,\chi\in{\cal D}$ 
$$\lim_n n\mathbb{E}[(\varphi(X)-\varphi(X_n))\chi(X_n)]=
\lim_n n\mathbb{E}[(X-X_n)\varphi^\prime(X_n)\chi(X_n)+\frac{1}{2}(X-X_n)^2\varphi^{\prime\prime}(X_n)\chi(X_n)]$$
and 
$$
\begin{array}{rl}
\mathbb{E}(X-X_n)^2]&={\displaystyle\mathbb{E}\sum_{k=n}^\infty\frac{\mathbb{E}[(1_{U_{k+1}\leq X_k}-X_k)^2|{\cal F}_k]}{(k+3)^2}}\\
\\
&=\sum_{k=n}^\infty\frac{1}{(k+3)^2}\mathbb{E}[X_k(1-X_k)]\\
&\sim \frac{1}{6n}
\end{array}
$$
because $X_n\rightarrow X$ and $X$ is uniformly distributed on $[0,1]$, as easily verified.

We obtain 
$$\lim_n n\mathbb{E}[(\varphi(X)-\varphi(X_n))\chi(X_n)]=\frac{1}{12}\mathbb{E}[\varphi^{\prime\prime}(X)\chi(X)]$$
$$\lim_n n\mathbb{E}[(\varphi(X)-\varphi(X_n)^2]=\frac{1}{6}\mathbb{E}[\varphi^{\prime 2}(X)].$$
Hence (H1) to (H3) are fulfilled
$
\underline{A}[\varphi]=\frac{1}{12}\varphi^{\prime\prime}$ and
$\widetilde{A}[\varphi]=\frac{1}{12}\varphi^{\prime\prime}
$
so that $\underline{A}=\overline{A}$ and $\Abar=0$. The limit error structure is the uniform error structure on $[0,1]$. This analysis 
could be easily extended to any configuration of Polya's urn, mutatis mutandis.

More generally, this kind of asymptotic behavior appears, under  regularity assumptions, for the approximation between a martingale and its limit.

Let $M_n=\sum_{i=1}^n X_i$ be a martingale w.r. to the filtration ${\cal F}_n$. Let us suppose $M_n$ centered, square 
integrable s.t. $\sum_{i=1}^\infty\mathbb{E} X_i^2<+\infty$ and let us put $\sigma_n^2=\sum_{i=n+1}^\infty\mathbb{E} X_i^2$.\\

\noindent{\bf Proposition  9. }{\it Supposing $\frac{1}{\sigma_n}\sup_{i> n}|X_i|\rightarrow 0$ in probability, 
$\mathbb{E}[\frac{1}{\sigma_n^2}\sup_{i>n}X_i^2]$ bounded in $n$, $\mathbb{E}[(\frac{1}{\sigma_n^2}(\sum_{i>n}X_i)^2)^p]$ bounded in $n$ for some $p>1$,
and $\frac{1}{\sigma_n^2}\sum_{i>n}X_i^2\rightarrow \zeta^2$ in probability, then

a) $\frac{1}{\sigma_n}\sum_{i>n}X_i \stackrel{d}{\Longrightarrow}Z$, where $Z$ has for characteristic 
function $\mathbb{E}[e^{-\frac{1}{2}\zeta^2t^2}]$

b) hypothesis } (H2) {\it is satisfied and $\forall\varphi,\chi\in{\cal C}^\infty_K$,
 
$\lim_n\frac{1}{\sigma_n^2}\mathbb{E}[(\varphi(M_\infty)-\varphi(M_n))\chi(M_n)]
=\frac{1}{2}\mathbb{E} Z^2\mathbb{E}[\varphi^{\prime\prime}(M_\infty)\chi(M_\infty)]$ what gives

$\underline{A}[\varphi]=\frac{1}{2}\mathbb{E}[Z^2]\,\varphi^{\prime\prime}.$

c) $\lim_n\frac{1}{\sigma_n^2}\mathbb{E}[(\varphi(M_\infty)-\varphi(M_n))^2]=\mathbb{E}[Z^2]\mathbb{E}[\varphi^{\prime 2}(M_\infty)]$, hypothesis } (H3)
{\it is satisfied as soon as the law of $M_\infty$ satisfies the Hamza condition, then $\underline{A}=\overline{A}=\widetilde{A}$ and 
$\Abar=0$.}\\

\noindent{\bf Proof.} The proposition is a direct consequence of a result of Hall and Heyde ([14] \S 3.5 p.76 {\it et seq.}).\\

\noindent{\bf II.3. Conditionally Gaussian case}\\

Let us begin with the finite dimensional case before applying the approach to processes.\\

{\bf II.3.a. } Let $Y$ be a r.v. with values in $\mathbb{R}^d$, $V=(V_{ij})_{i,j=1,\ldots,d}$ be an application from $\mathbb{R}^d$ into symmetric positive $d\!\times\!d$-matrices,
$\xi_j$ be r.v. with values in $\mathbb{R}^d$ which conditionally given $Y=y$ are i.i.d. Gaussian with common law ${\cal N}_d(y,V(y))$. We consider 
$$Y_n=\frac{1}{n}\sum_{j=1}^n \xi_j$$ as approximation of $Y$. We take ${\cal D}={\cal L}\{x\mapsto e^{i<u,x>},\;u\in\mathbb{R}^d\}$ and 
$\alpha_n=n$.\\

\noindent{\bf Lemma 1. }{\it Let us suppose $\mathbb{E}[{\mbox{\rm trace}}(V(Y))]<+\infty$, then

1) $\lim_n n\mathbb{E}[(\varphi(Y_n)-\varphi(Y))^2]
=\mathbb{E}_Y[\sum_{i,j=1}^d\varphi^\prime_i\varphi^\prime_j V_{ij}]\quad\forall\varphi\in{\cal D}$,

2) the hypothesis }(H1) {\it is satisfied and $\overline{A}[\varphi](y)
=\frac{1}{2}\sum_{i,j=1}^d V_{ij}(y)\varphi^{\prime\prime}_{ij}(y)\quad\forall\varphi\in{\cal D}$.}\\

The proof proceeds without difficulties as in the preceding example by finite expansions of the exponential functions.\\

\noindent{\bf Lemma 2. }{\it Let us suppose  $\mathbb{E}[{\mbox{\rm trace}}(V(Y))]<+\infty$ and the following regularity condition: 
$\forall i,j$, the measure $V_{ij}\mathbb{P}_Y$ has a partial derivative $\partial_j(V_{ij}Ñ\mathbb{P}_Y)$ in the sense of distributions which is a bounded measure
absolutely continuous w.r. to $\mathbb{P}_Y$, say $\rho_{ij}\mathbb{P}_Y$, then as soon as $V_{ij}$ and $\rho_{ij}\in L^2(\mathbb{P}_Y)$, } (H1) {\it to} (H3) {\it are fulfilled and }
$$
\begin{array}{l}
\underline{A}[\varphi]=\frac{1}{2}\sum_{ij}V_{ij}\varphi^{\prime\prime}_{ij}+\sum_{ij}\rho_{ij}\varphi^\prime_j\\
\widetilde{A}[\varphi]=\frac{1}{2}\sum_{ij}V_{ij}\varphi^{\prime\prime}_{ij}+\frac{1}{2}\sum_{ij}\rho_{ij}\varphi^\prime_j\\
\Abar[\varphi]=-\frac{1}{2}\sum_{ij}\rho_{ij}\varphi^\prime_j
\end{array}
$$
{\bf Proof.} The condition of the statement allows to perform an integration by parts in the limit obtained in lemma 2. 
That gives (H3) hence (H2) as well.\hfill$\diamond$\\

\noindent{\bf Remark 9. } There are several sufficient conditions in order that the form $\hat{\cal E}[u,v]=\mathbb{E}_Y[\sum_{ij}u^\prime_i v^\prime_j V_{ij}]$ be closable
on ${\cal C}^\infty_K(\mathbb{R}^d)$ (cf. [13] chapter 3 \S 3.1 and [23] chapter II \S 2). Suppose such a condition holds, then by the argument
of remark 2 the hypothesis (H3) is equivalent to ${\cal D}\subset {\cal DA}$ where $({\cal A,DA})$ is 
the generator of the smallest closed extension of $(\hat{\cal E},{\cal C}^\infty_K(\mathbb{R}^d))$.\\

{\bf II.3.b. }Let us apply this to the approximation of processes. Let $Y$ be a real process indexed by a set $T$. Let us consider $Y$ as a 
measurable map from $(\Omega,{\cal A},\mathbb{P})$ into $(E,{\cal F})=(\mathbb{R}^T,({\cal B}(\mathbb{R}))^{\otimes T})$. Let $y_t$ be the coordinate mappings from
$E$ on $\mathbb{R}$. We consider the algebra
$${\cal D}={\cal L}\{e^{iu_1y_{t_1}+\cdots+iu_ky_{t_k}},u_j\in\mathbb{R},t_j\in T\}.$$
Thanks to the monotone class theorem, ${\cal D}$ is a dense algebra in $L^2(E,{\cal F},\mathbb{P}_Y)$. We put $Y_t=y_t\circ Y$.

Let $\xi^j=(\xi^j_t)_{t\in T}$ be a sequence of real processes such that, conditionally given $Y$ the $\xi^j$ are independent with the same Gaussian law with
$$\mathbb{E}[\xi^1_t|Y]=Y_t$$
$$\mathbb{E}[(\xi^1_s-Y_s)(\xi^1_t-Y_t)|Y]=C_{s,t}(Y_s,Y_t)$$
the function $C_{s,t}(x_1,x_2)$ and the process $Y$ being such that $\mathbb{E}[C_{s,t}(Y_s,Y_t)]<+\infty\quad\forall s,t\in T$.

We approximate $Y$ by the process $Y_n=\frac{1}{n}\sum_{j=1}^n\xi^j$. The results obtained in the finite dimensional case give the following proposition :\\

\noindent{\bf Proposition 10. } {\it If the marginal laws of $Y$ $\mathbb{P}_{(Y_{t_1},\ldots,Y_{t_k})}(dy_{t_1}\ldots dy_{t_k})$ possess partial derivatives in the sense
of ditributions $\frac{\partial}{\partial y_{t_i}}$ which are bounded measures absolutely continuous w.r. to $\mathbb{P}_{(Y_{t_1},\ldots,Y_{t_k})}$, say 
$\lambda_{ij}(y_{t_1},\ldots,y_{t_k})\mathbb{P}_Y$, then hypotheses } (H1) {\it to} (H3) {\it are verified and for $\varphi\in{\cal D}$ we have}
$$
\begin{array}{l}
\overline{A}[\varphi](y_{t_1},\ldots,y_{t_k})=\frac{1}{2}\sum_{i,j=1}^k C_{t_i,t_j}(y_{t_i},y_{t_j})\frac{\partial^2 \varphi}{\partial y_{t_i}\partial y_{t_j}}\\
\widetilde{A}[\varphi](y_{t_1},\ldots,y_{t_k})=\frac{1}{2}\sum_{i,j=1}^k C_{t_i,t_j}(y_{t_i},y_{t_j})\frac{\partial^2 \varphi}{\partial y_{t_i}\partial y_{t_j}}
+\sum_{i,j}\lambda_{ij}(y_{t_1},\ldots,y_{t_k})\frac{\partial \varphi}{\partial y_{t_j}}\\
\Gamma[\varphi]=\sum_{i,j=1}^k C_{t_i,t_j}\frac{\partial \varphi}{\partial y_{t_i}}\frac{\partial \varphi}{\partial y_{t_j}}.
\end{array}
$$
Under the hypotheses of proposition 10, the form $({\cal E},\mathbb{D})$ is local and theorem 2 on asymptotic error calculus applies. Let us also remark
that if we define the  operator $\#$ on ${\cal D}$ by
$$\varphi^\#=\sum_{j=1}^k\frac{\partial \varphi}{\partial y_{t_j}}\circ Y.(\xi^1_{t_j}-Y_{t_j})$$
when $\varphi$ depends only on $y_{t_1},\ldots,y_{t_k}$, we have
$$\mathbb{E}[(\varphi^\#)^2|Y]=\Gamma[\varphi](Y)$$ and $\int(\varphi^\#)^2d\mathbb{P}_Y(y)=2{\cal E}[\varphi]$. It follows that the operator $\#$ extends
uniquely to $\mathbb{D}$ in a closed operator satisfying for $F\in{\cal C}^1\cap Lip(\mathbb{R}^p,\mathbb{R})$
$$(F(f_1,\ldots,f_p))^\#=\sum_{j=1}^pF^\prime_j(f_1,\ldots,f_p)f^\#_j\qquad\forall f_1,\ldots,f_p\in\mathbb{D}$$
In other words $\#$ plays the role of a gradient w.r. to the Dirichlet form $({\cal E},\mathbb{D})$.\\

\noindent{\bf Special case 1.} $Y$ is a real process indexed by $\mathbb{R}_+$, and $\xi^j_t=Y_t+W^j_t$ where $W^j_t$ are 
independent standard Brownian motions independent of $Y$. If the marginal laws of $Y$ have densities $f_{t_1,\ldots,t_k}$ s.t. 
$\frac{\partial f_{t_1,\ldots,t_k}}{\partial y_{t_i}}=\lambda_i(y_{t_1},\ldots,y_{t_k})f_{t_1,\ldots,t_k}$ with $\lambda_i\in L^1(\mathbb{R}^k)$, the construction
applies and yields a Dirichlet form with square field operator $\Gamma$ s.t. 
$$\Gamma[\varphi](Y)=-\sum_{ij} u_iu_j \,t_i\wedge t_j \varphi^2(Y) \;\mbox{ for }\; \varphi(Y)=e^{i\sum_ku_kY_{t_k}}.$$

Suppose $Y$ possesses second order moments, then the linear forms 
$\ell=\sum_{p=1}^ka_p(Y_{t_{p+1}}-Y_{t_p})$ are in $\mathbb{D}$ and $\Gamma[\ell]
=\sum_{p=1}^k a_p^2(t_{p+1}-t_p)$, so that on step functions $f$ 
$$\Gamma[\int f\,dy]=\int f^2(s)\,ds.$$
Thus this error structure may be called {\it the Ornstein-Uhlenbeck structure on the process $Y$}.

\noindent{\bf Special case 2.} With the above notation, let us assume that 
$$\xi^j_t=Y_t-h(t)W^j_t$$ where $h$ is a deterministic function. With the same hypotheses as above, we have on step functions $f$
$$\Gamma[\int f\,dy]=\mathbb{E}[(\int f(t)d(h(t)W_t))^2]$$

\noindent{\bf Special case 3.} Suppose eventually 
$$\xi^j_t=Y_t+\int_0^th(s)\,dW_s^j$$
then $C_{t_i,t_j}=\int_0^{t_i\wedge t_j}h^2(s)\,ds$ and $\Gamma[\int f \,dy]=\int f^2(s)h^2(s)\,ds.$
We obtain a structure which may be called {\it the Ornstein-Uhlenbeck structure with weight $h$ on the process $Y$.}\\

\noindent{\bf II.4. Natural inaccuracy of the Brownian motion simulated by the Donsker theorem.}\\

We begin with the simplest case of one dimensional marginal laws which is here nothing else than the central limit theorem.\\

{\bf II.4.a. Natural inaccuracy in the central limit theorem.}

Let be $S_p=\sum_{i=1}^pV_i$ where the randon variables $V_i$ are i.i.d. centered with variance $\sigma^2$. We consider two indices
$m$ and $n$ linked by the relation
$$n=n(m)=m+k(m)\quad{\mbox{with }}\quad \theta\sqrt{m}\leq k(m)\leq\frac{1}{\theta}\sqrt{m}$$
for a $\theta\in]0,1[$.

Let us observe some evident properties : $\lim_{m\rightarrow\infty}\frac{m}{n}=1$; $\frac{n-m}{\sqrt{m}}\in[\theta,\frac{1}{\theta}]$; 
$\sqrt{n}\geq\sqrt{m}\sqrt{1+\theta}\geq\theta\sqrt{1+\theta}(n-m)$; $\frac{1}{\sqrt{m}}-\frac{1}{\sqrt{n}}\leq \frac{1}{2\theta m}$.

Writing $\frac{1}{\sqrt{n}}S_n-\frac{1}{\sqrt{m}}S_m=\frac{1}{\sqrt{n}}(S_n-S_m)+(\frac{1}{\sqrt{n}}-\frac{1}{\sqrt{m}})S_m$ and using these properties, 
shows that $\frac{1}{\sqrt{n}}S_n-\frac{1}{\sqrt{m}}S_m\rightarrow 0$ a.s.

We consider the mutual approximation of $\frac{1}{\sqrt{m}}S_m$ and $\frac{1}{\sqrt{n}}S_n$ (which is an obvious extension of the
 framework of part I). That is for $\overline{A}$ we study
\begin{equation}
\alpha(m)\mathbb{E}[(\varphi(\frac{1}{\sqrt{m}}S_m)-\varphi(\frac{1}{\sqrt{n}}S_n))\chi(\frac{1}{\sqrt{n}}S_n)]
\end{equation}
and for $\underline{A}$ we study
\begin{equation}
\alpha(m)\mathbb{E}[(\varphi(\frac{1}{\sqrt{n}}S_n)-\varphi(\frac{1}{\sqrt{m}}S_m))\chi(\frac{1}{\sqrt{m}}S_m)]
\end{equation}
with $\alpha(m)=\frac{m}{k(m)}$ (so that $\theta\sqrt{m}\leq\alpha(m)\leq \frac{1}{\theta}\sqrt{m}$). For the algebra ${\cal D}$ we take the linear combinations of imaginary exponentials.\\

\noindent{\bf Proposition 11. }{\it Suppose the $V_i$'s possess a third order moment, then hypotheses} (H1) {\it to} (H3) {\it are fulfilled and for 
$\varphi\in{\cal D}$
$$\overline{A}[\varphi](x)=\underline{A}[\varphi](x)=\widetilde{A}[\varphi](x)=\frac{1}{2}\sigma^2 \varphi^{\prime\prime}
-\frac{1}{2}x\varphi^\prime.$$
The Dirichlet form is the Ornstein-Uhlenbeck form on $\mathbb{R}$ (endowed with the normal law ${\cal N}(0,\sigma^2)$).}\\

\noindent{\bf Proof.} For $\overline{A}$, taking $\varphi(x)=e^{iux}$, $\chi(x)=e^{ivx}$ we have to look at 
$$J_m=\alpha(m)\mathbb{E}[e^{i(\frac{u}{\sqrt{m}}+\frac{v}{\sqrt{m}})S_m+\frac{iv}{\sqrt{n}}(S_n-S_m)}-e^{i\frac{u+v}{\sqrt{n}}S_n}].$$
Let $\xi(t)=\mathbb{E}[e^{itV_1}]$ be the characteristic function of the $V_i$'s
$$J_m=\alpha(m)[\xi(\frac{u}{\sqrt{m}}+\frac{v}{\sqrt{n}})^m(\xi(\frac{v}{\sqrt{n}}))^{n-m}-(\xi(\frac{u+v}{\sqrt{n}}))^n]$$
using the existence of a third moment we have
\begin{equation}
\log \xi(t)=-\frac{t^2}{2}\sigma^2(1+at+t\,o(1))
\end{equation}
and this allows to write 
$J_m=\alpha(m)[e^{-\frac{\sigma^2}{2}K_m}-e^{-\frac{\sigma^2}{2}L_m}]$ with

$$K_m=(u+v)^2+\frac{a}{\sqrt{m}}(u+v)^3+2uv(\frac{\sqrt{m}}{\sqrt{n}}-1)+\frac{1}{\sqrt{m}}\,o(1)$$ 

\noindent and  $L_m=(u+v)^2+\frac{a}{\sqrt{m}}(u+v)^3+\frac{1}{\sqrt{m}}\,o(1)$. 
This yields 

$$J_m=\alpha(m)e^{-\frac{\sigma^2}{2}(u+v)^2}[(-\frac{\sigma^2}{2})(2uv(\sqrt{\frac{m}{n}}-1)+\frac{1}{\sqrt{m}}\,o(1))]$$

\noindent hence 

$$\lim_m J_m=\frac{1}{2}\sigma^2uve^{-\frac{\sigma^2}{2}(u+v)^2}
=-\frac{\sigma^2}{2}\int_{\mathbb{R}}\varphi^\prime\chi^\prime\,d{\cal N}(0,\sigma^2)$$

\noindent what gives easily the proposition.\hfill$\diamond$\\

\noindent{\bf II.4.b. The Donsker case. }

Let the $V_i$'s be as before and 
$$X_n(t)=\frac{1}{\sqrt{n}}\left(\sum_{k=1}^{[nt]}V_k+(nt-[nt]V_{[nt]+1}\right)$$
for $t\in[0,1]$, $[nt]$ denoting the entire part of $nt$.

The laws of the variables $X_n$ are probability measures on ${\cal C}([0,1])$ as well as their limit in law which is a Brownian motion by Donsker theorem.

For the algebra ${\cal D}$ we take the linear combinations of exponential of the form 
$\varphi(X)=\exp{\{iX(f)\}}$ where $X(f)=\int_0^1 f(s)\,dX(s)$ and with $f\in{\cal C}^1$ in order that 
$\int_0^1 f(s)\,dX(s)$ may be defined as 
$X(1)f(1)-X(0)f(0)-\int_0^1X(s)df(s)$ for the general coordinate process $X(s)$ on ${\cal C}[0,1]$.
As easily seen the algebra ${\cal D}$ is dense in $L^2({\cal C}([0,1]),\mu)$ $\mu$ being the Wiener measure.

Thus we have $X_n(f)=\sqrt{n}\sum_{k=0}^{n-1}\int_{k/n}^{(k+1)/n}f(s)\,ds\,V_{k+1}$ and for studying the operator $\overline{A}$ we have
 to look at
$$M_m=\alpha(m)\mathbb{E}[(\varphi(X_m)-\varphi(X_n))\chi(X_n)]=\alpha(m)\mathbb{E}[(e^{iX_m(f)}-e^{iX_n(f)})e^{iX_n(g)}]$$
$$=\alpha(m)\mathbb{E}\left[\left(e^{i\sqrt{m}\sum_{k=0}^{m-1}\int_{\frac{k}{m}}^{\frac{k+1}{m}}f(s)dsV_{k+1}}-
e^{i\sqrt{n}\sum_{k=0}^{n-1}\int_{\frac{k}{n}}^{\frac{k+1}{n}}f(s)dsV_{k+1}}\right)
e^{i\sqrt{n}\sum_{k=0}^{n-1}\int_{\frac{k}{n}}^{\frac{k+1}{n}}g(s)dsV_{k+1}}\right].$$
We take as before $\alpha(m)=m/k(m)$.\\

\noindent{\bf Proposition 12. } {\it Suppose the $V_i$'s possess a third order moment, then hypotheses } 
(H1) {\it to} (H3) {\it are fulfilled. We have $\overline{A}=\underline{A}=\widetilde{A}$ on ${\cal D}$. The Dirichlet form is the Ornstein-Uhlenbeck
form on the Wiener space (with a Brownian motion s.t. $<B>_t=\sigma^2 t$) normalized so that the square field operator satisfies}
$\Gamma[\int_0^1 h(s)\,dB_s]=\int_0^1h^2(s)\,\sigma^2\,ds\quad\forall h\in L^2([0,1]).$\\

\indent Since the Dirichlet form is local, some limits are automatically obtained (theorem 2). Since $\Abar=0$, the theoretical and practical variances
coincide (prop. 1).

\noindent{\bf Proof.} For studying $\overline{A}$ we consider the quantity $M_m$ defined above. By the third moment assumption,
the characteristic function $\xi(t)$ of the $V_i$'s satisfies (7) and we can write
$$M_m=\alpha(m)\left[\prod_{j=m}^{n-1}\xi(\sqrt{n}\int_{\frac{j}{n}}^{\frac{j+1}{n}}\!g)
\prod_{k=0}^{m-1}\xi(\sqrt{m}\int_{\frac{k}{m}}^{\frac{k+1}{m}}\!f+\sqrt{n}\int_{\frac{k}{n}}^{\frac{k+1}{n}}\!g)
-\prod_{j=m}^{n-1}\xi(\sqrt{n}\int_{\frac{j}{n}}^{\frac{j+1}{n}}\!(f+g)\right]$$
\begin{equation}
=\alpha(m)\left[e^{-\frac{\sigma^2}{2}N_m}-e^{-\frac{\sigma^2}{2}P_m}\right]
\end{equation}
with
$$N_m=\sum_{j=m}^{n-1}n(\int_{\frac{j}{n}}^{\frac{j+1}{n}}\!g)^2+an\sqrt{n}(\int_{\frac{j}{n}}^{\frac{j+1}{n}}\!g)^3(1+o(1))\qquad\qquad\qquad\qquad\qquad\qquad\qquad$$
$$\qquad\qquad\qquad+\sum_{k=0}^{m-1}(\sqrt{m}\int_{\frac{k}{m}}^{\frac{k+1}{m}}\!f+\sqrt{n}\int_{\frac{k}{n}}^{\frac{k+1}{n}}\!g)^2
+a(\sqrt{m}\int_{\frac{k}{m}}^{\frac{k+1}{m}}\!f+\sqrt{n}\int_{\frac{k}{n}}^{\frac{k+1}{n}}\!g)^3(1+o(1))$$
$$P_m=\sum_{k=0}^{n-1}n(\int_{\frac{j}{n}}^{\frac{j+1}{n}}\!(f+g))^2+an\sqrt{n}(\int_{\frac{j}{n}}^{\frac{j+1}{n}}\!(f+g))^3(1+o(1)).\qquad\qquad\qquad\qquad\qquad\qquad$$
Using $m\sum_{j=0}^{m-1}(\int_{\frac{j}{n}}^{\frac{j+1}{n}}\!f)^2=\int_0^1f^2+\frac{1}{m}O(1)$ and 
$\sup_j|\sqrt{n}\int_{\frac{j}{n}}^{\frac{j+1}{n}}g|\leq\frac{1}{\sqrt{n}}\|g\|_\infty$ we obtain
$$N_m=\int_0^1g^2+\int_0^1f^2+2\sum_{k=0}^{m-1}\sqrt{m}\sqrt{n}\int_{\frac{k}{m}}^{\frac{k+1}{m}}\!f
\int_{\frac{k}{n}}^{\frac{k+1}{n}}\!g+\frac{a}{\sqrt{m}}\int_0^1(f+g)^3(1+o(1))+\frac{1}{n}O(1)$$
and
$$P_m=\int_0^1(f+g)^2+\frac{a}{\sqrt{m}}\int_0^1(f+g)^3(1+o(1))+\frac{1}{n}O(1)$$
Putting these expressions in (8) leads to 
$$M_m=\alpha(m)\left(\exp[-\frac{\sigma^2}{2}\int_0^1(f+g)^2]\right)\left[(-\sigma^2)(\sqrt{\frac{m}{n}}-1)\int_0^1fg+\frac{1}{\sqrt{m}}o(1)\right]$$
Eventually, for $\varphi(X)=\exp[i\int_0^1fdX]$ and $\chi(X)=\exp[i\int_0^1gdX]$ we get 
$$\alpha(m)\mathbb{E}[(\varphi(X_m)-\varphi(X_n))\chi(X_n)]\rightarrow\frac{\sigma^2}{2}\exp[-\frac{\sigma^2}{2}\int_0^1(f+g)^2\,ds]\int_0^1fg\,ds.$$
In order to recognize the obtained limit, let $\Gamma_{ou}$ be the Ornstein-Uhlenbeck square field operator on the standard Wiener 
space (s.t. $\Gamma_{ou}[\int h\,dB]=\int h^2ds$). We have by the functional calculus
$$\mathbb{E}\Gamma_{ou}[e^{i\int_0^1fdB},e^{i\int_0^1 gdB}]=-\int fg\,ds\mathbb{E}[e^{i\int(f+g)dB}]=-\int fg\,ds\exp[-\frac{1}{2}\int(f+g)^2\,ds].$$
It follows that for a Wiener measure s.t. $<B>_t=\sigma^2 t$ and an Ornstein-Uhlenbeck structure s.t. $\Gamma[\int_0^1hdB]=\int h^2d<B>$ whose
generator will be denoted $(A,{\cal D}A)$, we have ${\cal D}\subset{\cal D}A$ and 

$$\begin{array}{rl}
<A\varphi,\chi>&=-{\cal E}[\varphi,\chi]=-\frac{1}{2}\mathbb{E}[\Gamma[\varphi,\chi]]\\
&=\frac{\sigma^2}{2}\int_0^1fg\,ds\exp[-\frac{\sigma^2}{2}\int_0^1(f+g)^2\,ds]=<\overline{A}\varphi,\chi>.
\end{array}
$$
The operator $\overline{A}$ is therefore symmetric on ${\cal D}$, which implies $\overline{A}=\underline{A}=\widetilde{A}$ 
and the proposition is proved.\hfill$\diamond$\\

\noindent{\bf Comment.} As noted already by Louis Bachelier, assets quoted on the organized markets {\it look like} Brownian paths. 
This is displayed in any course in mathematical finance in order to introduce modelling by diffusion processes and stochastic calculus. Nevertheless
some concrete facts prevent this ressemblance from being accurate at microscopic scale. First because the spot is only defined at discrete
instants, second because a Brownian path possesses ideal properties (like the fact that it cuts uncountably many times every level that it reaches) that
cannot be verified by material recordings.

So that in order to be completely pragmatic, we might replace in financial models any Brownian motion by a random walk by application of Donsker theorem
with $n$ sufficiently large and consider the Brownian motion of the model {\it is nothing else} than a class of such sufficiently fine 
random walks. An infinite precision for stochastic calculus in finance is therefore a priori absurd and, by the results of this section, we may represent
the intrinsic fuzzyness of these computations by the Ornstein-Uhlenbeck form on the Wiener space. This is a justification of 
the approach proposed in [4].\\

\noindent{\bf II.5. Empirical laws and natural inaccuracy of the Brownian bridge}\\

If $(V_n)$ are i.i.d. real random variables $0\leq V_n\leq 1$ with distribution function $F$ and if $F_n(x) =\frac{1}{n}\sum_{i=1}^n 1_{V_i\leq x}$ is the
empirical distribution function, then $\sqrt{n}(F_n-F)$ converges in law on the Skorohod space to a transformed Brownian bridge 
$B_{F(x)}-F(x)B_1$ (see for instance [3] p. 141).

Considering the $V_i$'s are simulated by the inversion method shows that this result is a consequence of the special case where the $V_i$'s are 
uniformly ditributed on $[0,1]$. From now on, we restrict to this case. Putting $Z_n(x)=\sqrt{n}(F_n(x)-x)$ we are interested in the limit 
$R_m=\alpha(m)\mathbb{E}[(\varphi(Z_m)-\varphi(Z_n))\chi(Z_n)]$ for $m$ and $n$ linked as in the preceding example, 
with also $\alpha(m)=m/k(m)$ and for $\varphi,\chi\in{\cal D}$ where ${\cal D}$ is the algebra of linear combinations 
of imaginary exponentials of the form $\varphi(Z)=\exp\{i\int_0^1f(s)dZ(s)\}=\exp\{-i\int_0^1Z(s)df(s)\}$ for $f\in{\cal C}^1([0,1])$. 

Thus
$$R_m=\alpha(m)\mathbb{E}\left[\left(e^{i\frac{1}{\sqrt{m}}\sum_{k=1}^m(f(V_k)-\int f)}
-e^{i\frac{1}{\sqrt{n}}\sum_{k=1}^n(f(V_k)-\int f)}\right)e^{i\frac{1}{\sqrt{n}}\sum_{k=1}^n(g(V_k)-\int g)}\right].$$
Putting $\tilde{f}=f-\int_0^1 fds$ and $\tilde{g}=g-\int_0^1 gds$ and denoting $\eta$ and $\zeta$ the characteristic functions of 
$\tilde{g}$ and $\tilde{f}+\tilde{g}$ gives
$$R_m=\alpha(m)\left[(\rho(m))^m\left((\eta(\frac{1}{\sqrt{n}}))^{n-m}-(\zeta(\frac{1}{\sqrt{n}}))^n\right)\right]$$
with $\rho(m)=\mathbb{E}[e^{i\frac{1}{\sqrt{m}}\tilde{f}+i\frac{1}{\sqrt{n}}\tilde{g}}]$.

The estimates
$$\begin{array}{l}
\rho(m)=1-\frac{1}{2m}\mathbb{E}[(\tilde{f}+\sqrt{\frac{m}{n}}\tilde{g})^2]
+(\frac{i}{\sqrt{m}})^3\frac{1}{6}\mathbb{E}[(\tilde{f}+\sqrt{\frac{m}{n}}\tilde{g})^3]+\frac{1}{m^2}O(1)\\
\zeta(\frac{1}{\sqrt{n}})=1-\frac{1}{2n}\mathbb{E}[(\tilde{f}+\tilde{g})^2]
+(\frac{i}{\sqrt{n}})^3\frac{1}{6}\mathbb{E}[(\tilde{f}+\tilde{g})^3]+\frac{1}{m^2}O(1)\\
\eta(\frac{1}{\sqrt{n}})=1-\frac{1}{2n}\mathbb{E}[\tilde{g}^2]+(\frac{i}{\sqrt{n}})^3\frac{1}{6}\mathbb{E}[\tilde{g}^3]+\frac{1}{m^2}O(1)
\end{array}
$$
allow to obtain
$$\begin{array}{rl}
\lim_m R_m&=e^{-\frac{1}{2}\mathbb{E}[(\tilde{f}+\tilde{g})^2]}(-\frac{1}{2}\mathbb{E}[\tilde{g}^2]+\frac{1}{2}cov(g,f+g))\\
&=\frac{1}{2}e^{-\frac{1}{2}\mathbb{E}[(\tilde{f}+\tilde{g})^2]}cov(f,g).
\end{array}
$$
In order to recognize the limit, let $\Gamma_{ou}$ be as before the Ornstein-Uhlenbeck square field operator on the Wiener space
and $Z_t=B_t-tB_1$. We have 
$$\Gamma_{ou}[e^{i\int fdZ},e^{i\int gdZ}]=-e^{i\int(f+g)dZ}\Gamma_{ou}[\int fdZ,\int gdZ]
=-e^{i\int(f+g)dZ}(\int fgds-\int fds\int gds)$$
and consequently
$$\begin{array}{rl}
{\cal E}_{ou}[\varphi,\chi]=\frac{1}{2}\mathbb{E}[\Gamma_{ou}[\varphi,\chi]]&
=-\frac{1}{2}e^{-\frac{1}{2}[\int(f+g)^2ds-(\int(f+g)ds)^2]}(\int fgds-\int fds\int gds)\\
&=\lim_m R_m.
\end{array}
$$
As before denoting $A$ the Ornstein-Uhlenbeck operator, we see that $\overline{A}=A$ on ${\cal D}$ and therefore $\overline{A}=\underline{A}=\widetilde{A}$, 
the Dirichlet form is the  image of the Ornstein-Uhlenbeck
form on the Wiener space.\\

\noindent{\bf II.6. Erroneous empirical laws and generalized Mehler type structures on the Brownian bridge}\\

We still consider a sequence $V=(V")_{i\in\mathbb{N}}$ of i.i.d. random variables uniformly distributed on $[0,1]$ and the 
empirical distribution function $$F_n(x)=\frac{1}{n}\sum_{i=1}^n1_{V"\leq x}$$ but the problem that we tackle is different. We suppose that the law
of $V^1$ is not perfectly known. We assume that there is a sequence of r.v. $U^1_m$ approximating $V^1$ and copies $(U^i_m,V^i)$ of 
$(U^1_m,V^1)$ such that the sequence $(U^i_m,V^i)_{i\in\mathbb{N}}$ be i.i.d. and we suppose that the law of $U^1_m$ has support in $[0,1]$ with distribution
function $F^m$. We define the emprirical distribution function
$$F^m_n(x)=\frac{1}{n}\sum_{j=1}^n1_{U^j_m\leq x}.$$
We are interested in the approximation 
$$Z^m_n=\frac{1}{\sqrt{n}}\sum_{j=1}^n(1_{U^j_m\leq x}-F^m(x))$$
of the process
$$Z_n=\frac{1}{\sqrt{n}}\sum_{j=1}^n(1_{V^i\leq x}-x).$$
We take ${\cal D}={\cal L}\{\varphi(z)=\exp[-i\int_0^1z(t)df(t)],\quad f\in{\cal C}^1([0,1])\}$ and we study
$$Q_{m,n}=\beta_m\mathbb{E}[(\varphi(Z^m_n)-\varphi(Z_n))^2]
=\beta_m\mathbb{E}[(\exp[i\frac{1}{\sqrt{n}}\sum_{i=1}^n\widetilde{f(U^i_m)}]-\exp[i\frac{1}{\sqrt{n}}\sum_{i=1}^n\widetilde{f(V^i)}])^2]$$
in which the symbol $\widetilde{\;}$ represents the centering operation.

Denoting $\eta$ [resp. $\eta_m$] the characteristic function of $\widetilde{f(V^1)}$ [resp. $\widetilde{f(U^i_m)}$], $\theta_m$ the characteristic
function of $\widetilde{f(U^1_m)}+\widetilde{f(V^1)}$, $\sigma^2$ [resp. $\sigma^2_m$] the variance of $f$ [resp. of 
$\widetilde{f(U^1_m)}$],
we have
$$Q_{m,n}=\beta_m[(\eta_m(\frac{2}{\sqrt{n}}))^n-2(\theta_m(\frac{1}{\sqrt{n}}))^n+(\eta(\frac{2}{\sqrt{n}}))^n]$$
The estimates
$$\begin{array}{l}
\eta_m(t)=1-\frac{t^2}{2}\sigma_m^2+t^2o(1)\\
\theta_m(t)=1-\frac{t^2}{2}\mathbb{E}[(\widetilde{f(U^1_m)}+\widetilde{f(V^1)})^2]+t^2o(1)\\
\eta(t)=1-\frac{t^2}{2}\sigma^2+t^2o(1)
\end{array}$$
give
$$Q_{m,n}=\beta_me^{-2\sigma^2}\left[e^{-2(\sigma_m-\sigma^2)+o(1)}-2e^{2\sigma^2-\frac{1}{2}\mathbb{E}[(\widetilde{f(U^1_m)}+\widetilde{f(V^1)})^2]+o(1)}
+\exp o(1)\right].$$
Remarking that $2\sigma^2-\frac{1}{2}\mathbb{E}[(\widetilde{f(U^1_m)}+\widetilde{f(V^1)})^2]=\frac{1}{2}\mathbb{E}[(\widetilde{f(U^1_m)}+\widetilde{f(V^1)})^2]
+\sigma^2-\sigma_m^2$, if we assume $\mathbb{E}[(\widetilde{f(U^1_m)}-\widetilde{f(V^1)})^2]\rightarrow_{m\rightarrow\infty} 0 $ and $\sigma^2-\sigma_m^2\rightarrow 0$,
we obtain
$$Q_{m,n}=\beta_me^{-2\sigma^2}(-\mathbb{E}[(\widetilde{f(U^1_m)}-\widetilde{f(V^1)})^2]+o(1))$$
we can state\\

\noindent{\bf Proposition 13. } {\it If there is a sequence $\beta_m\rightarrow +\infty$ s.t. 
$$\beta_m\mathbb{E}[(\widetilde{f(U^1_m)}-\widetilde{f(V^1)})^2]\rightarrow \frac{1}{2}e[f]$$ where
$e[.]$ is a quadratic form defined on ${\cal C}^1([0,1])$ closable in $L^2([0,1])$ (with a Dirichlet extension non necessarily local, cf. remark 2) and supposing
$\mathbb{E}[(\widetilde{f(U^1_m)})^2]\rightarrow \mathbb{E}[(\widetilde{f(V^1)})^2]$ then
$$\lim_{m,n\uparrow\infty}\beta_m\mathbb{E}[(e^{i\int fdZ^m_n}-e^{i\int fdZ_n})^2]=\frac{1}{2}e^{-2\sigma^2}e[f],$$
hypothesis} (H3) {\it is fulfilled and the limit Dirichlet form is the image by the Brownian bridge of the generalized Mehler
 type form on the Wiener space associated with the form $e[.]$
}(cf. [5] chapter VI \S 2.5 p113 et seq).\\

\noindent{\bf Proof.} The hypotheses of the proposition imply what we needed during the above computation. It suffices therefore, as before, to recognize
the limit as a closable form. But that comes from the functional calculus and the fact that the generalized Mehler type structure 
associated with the form $e[.]$ satisfies $\Gamma[\int_0^1 f(s)\,dB_s]=e[f]$.\hfill$\diamond$

For example if $U^1_m=\frac{1}{m}\sum_{p=1}^m(V^1+a(V^1)G_p)$ where the $G_p$ are i.i.d. reduced normal variables 
independent of $V^1$ and where $a$ is continuous, then by the lemma 1  of the conditionally Gaussian case and the Hamza condition, the form
$$e[f]=\lim_n n\mathbb{E}[(\widetilde{f(U^1_m)}-\tilde{f}(V^1))^2]=\int_0^1 f^{\prime 2}(x)\,a^2(x)\,dx$$ is closable and the proposition applies.
This generalized Mehler type structure satisfies
\begin{equation}
\Gamma[\int_0^1 h\,dB]=\int_0^1f^{\prime 2}(s)\,ds\quad\forall f\in H^1([0,1])
\end{equation}
this structure may be constructed elementarily as in [5] or by the white noise theory.\\

\noindent{\bf Remark 10. } There exists an extension of Donsker theorem to the case where the variables $V_k$ (notation of section II.4.b)
are erroneous with a functional convergence in the sense of Dirichlet forms (see [6] and [10]). The limit structure obtained is the Ornstein-Uhlenbeck 
structure. This is related to the fact that the perturbation involved in this 
approach is a
{\it transversal} perturbation of the random walk hence at the limit a transversal perturbation of the Brownian path
 (we will display this result in terms of an approximation procedure in section II.7 below).

Here instead, we change the law of the starting random variables : In the expression of $Z^m_n$
$$Z^m_n=\frac{1}{\sqrt{n}}\sum_{j=1}^n(1_{U^j_m\leq x}-F^m(x))$$ if we simulate $U^j_m$ by the inversion method we see that $Z^m_n$ has same law
as
$$\frac{1}{\sqrt{n}}\sum_{j=1}^n(1_{(F^m)^{-1}(X^j)\leq x}-F^m(x))=\frac{1}{\sqrt{n}}\sum_{j=1}^n(1_{X^j\leq F^m(x)}-F^m(x))$$
where  $(X^j)_{j\in\mathbb{N}}$ is a copy of  $(V^j)_{j\in\mathbb{N}}$. We see that when $F^m$ changes, the path is {\it longitudinally} perturbed
and so is the limit Brownian bridge. This explains a formula like (9).\\

\noindent{\bf II.7. Erroneous random walk and Donsker theorem.}\\

This example displays many similarities with examples II.4, II.5 and II.6. We give only the framework and the results.

Let $U^1$ be a centered square integrable r.v. approximated by $U^1_m$ also centered and square integrable. We suppose 
$$\alpha_m\mathbb{E}[(U_m-U)^2]\rightarrow\lambda.$$
Considering i.i.d. copies $(U^i_m,U^i)$ of $(U^1_m,U^1)$ we look at 
$$\begin{array}{l}
{\displaystyle X^n_m(t)=\frac{1}{\sqrt{n}}(\sum_{i=1}^{[nt]}U^i_m+(nt-[nt])U^{[nt]+1}_m)}\\
{\displaystyle X^n(t)=\frac{1}{\sqrt{n}}(\sum_{i=1}^{[nt]}U^i+(nt-[nt])U^{[nt]+1})}
\end{array}
$$
and we study $T_{m,n}=\alpha_m\mathbb{E}[(\varphi(X^n_m)-\varphi(X^n))^2]$ for $\varphi$ belonging to

$${\cal D}={\cal L}\{\exp[i\int_0^1f(s)\,dX_s], \quad f\in{\cal C}^\infty\mbox{ with support in }]0,1[\}.$$
Putting $\sigma^2=var(U^1)$ we find that 
$$\lim_{m,n\uparrow\infty}T_{m,n}=-\lambda(\int f^2\,ds) e^{-2\sigma^2\int f^2\,ds}$$
The limit Dirichlet form is once more the Ornstein-Uhlenbeck form on the Wiener space s.t. $<B>_t=\sigma^2 t$
and $\Gamma[\int_0^1 f\,dB]=\lambda\int_0^1 f^2(s)\,d<B>_s$.\\

\noindent{\bf Comment. } The Dirichlet-version of the Donsker theorem proved in [6] supposes 
the r.v. $U^1$ has a regular law allowing to carry a (non zero) Dirichlet form. This excludes the case of a discrete law. Here instead, we do not need
such a restriction and the present construction applies for instance to the Cox-Ross-Rubinstein model approximating the Black-Scholes model. But
the convergence here is weaker than the one used in [6].\\

\noindent{\bf II.8. Approximation of the Brownian motion defined through the Wiener integral, i.e. as centered orthogonal measure.}\\

Let $X^1_n$ be a centered square integrable real r.v. approximating the variable $X^1$ which is reduced Gaussian. Let $(X^i_n,X^i)_{i\in\mathbb{N}^\ast}$
be i.i.d. copies of $(X_n^1,X^1)$.

We assume $X^1_n$ and $X^1$ satisfy the hypothesis (H3) with $\alpha_n$ and an algebra ${\cal D}_0$ dense in $L^2({\cal N}(0,1))$ 
of bounded functions containing the constants and the imaginary exponentials: $\forall\varphi,\chi\in{\cal D}_0$ 
\begin{equation}
\alpha_n\mathbb{E}[(\varphi(X^1_n)-\varphi(X^1))(\chi(X^1_n)-\chi(X^1))]\rightarrow \lambda\int_{\mathbb{R}}\varphi^\prime\chi^\prime\,d{\cal N}(0,1).
\end{equation}
Let $\xi_k$ be an orthonormal basis of $L^2(E_1,{\cal F}_1,\mu_1)$ where $(E_1,{\cal F}_1,\mu_1)$ is a $\sigma$-finite measured space and 
let us consider the mappings
$$\begin{array}{l}
{\displaystyle f\in L^2(E_1,{\cal F}_1,\mu_1)\quad\stackrel{J_n}{\mapsto}\quad\sum_{k=1}^\infty<f,\xi_k>X^k_n}\\
{\displaystyle f\in L^2(E_1,{\cal F}_1,\mu_1)\quad\stackrel{J}{\mapsto}\quad\sum_{k=1}^\infty<f,\xi_k>X^k.}
\end{array}
$$
We consider $J_n$ as an approximation of $J$ and for the algebra ${\cal D}$ we choose
$$
\begin{array}{l}
{\cal D}={\cal L}\{\Phi\,:\,\theta\in\mbox{CORM}\mapsto\Phi(\theta)=e^{i<f,\theta>}\\
\hspace{3cm}\mbox{ where }f\in L^2(E_1,{\cal F}_1,\mu_1)
\mbox{ has a finite expansion on the basis }(\xi_k)\}
\end{array}
$$
here CORM denotes the set of centered orthogonal random measures on $L^2(E_1,{\cal F}_1,\mu_1)$.

We study the limit of $\Delta_n=\alpha_n\mathbb{E}[\Phi(J_n)-\Phi(J))^2]$ for $\Phi(\theta)=e^{i<f,\theta>}$ 
with $f=\sum_{q=1}^Q f_q\xi_d$.

We have
$$\begin{array}{rl}
\Delta_n&=\alpha_n\left(\prod_{k=1}^Q\mathbb{E} e^{2if_kX^1_n}-2\prod_{k=1}^Q\mathbb{E} e^{if_k(X^1_n+X^1)}+\prod_{k=1}^Q\mathbb{E} e^{2if_kX^1}\right)\\
&=\alpha_n\left(\prod_{k=1}^Q\mathbb{E} e^{2if_kX^1_n}-\prod_{k=1}^Q\mathbb{E} e^{if_k(X^1_n+X^1)}+\prod_{k=1}^Q\mathbb{E} e^{2if_kX^1}-\prod_{k=1}^Q\mathbb{E} e^{if_k(X^1_n+X^1)}\right)\\
&=\alpha_n(A_n+B_n)
\end{array}
$$
and we may write
$$A_n=\sum_{q=1}^Q\prod_{k=1}^{q-1}\mathbb{E} e^{2if_kX^1_n}\left(\mathbb{E} e^{2if_qX^1_n}-\mathbb{E} e^{2if_q(X^1_n+X^1)}\right)\prod_{k=q+1}^Q\mathbb{E} e^{2if_k(X^1_n+X^1)}$$
because the intermediate terms cancel and remain only the first and the last ones.

Let us assume now in addition that the pair $(X^1_n,X^1)$ converges in law to $(X^1,X^1)$. Then
$$
\begin{array}{rl}
\lim_n \alpha_n(A_n+B_n)&{\displaystyle=\sum_{q=1}^Q\prod_{k=1}^{q-1}\mathbb{E} e^{2if_kX^1}\lim_n\mathbb{E}[(e^{if_dX^1_n}-e^{if_qX^1})^2]\prod_{k=q+1}^Q\mathbb{E} e^{2if_kX^1}}\\
&{\displaystyle=\sum_{q=1}^Q\mathbb{E}[e^{2i\sum_{k=1}^Qf_kX^k}]\lim_n\alpha_n\frac{\mathbb{E}[(e^{if_qX^1_n}-e^{if_qX^1})^2]}{\mathbb{E}[e^{2if_qX^1}]}}
\end{array}
$$
and by the assumption (10) this is nothing else than 
$$\lim_n\Delta_n=\mathbb{E}[e^{2i\sum_{k=1}^Qf_kX^k}]\sum_{q=1}^Q(-\lambda f_q^2)=\mathbb{E}[e^{2iJ(f)}]\|f\|^2$$
We recognize once more the Ornstein-Uhlenbeck structure on the abstract Wiener space defined by $J$. We can state\\

\noindent{\bf Proposition 14. } {\it If $X^1_n$ is an approximation of $X^1$ satisfying }(10) {\it and if $(X^1_n,X^1)\stackrel{d}{\Rightarrow}(X^1,X^1)$
then the approximation $J_n$ of the centered orthogonal random measure $J$ satisfies } (H3) {\it on ${\cal D}$ and yields the Ornstein-Uhlenbeck form.}\\

We would prove easily following the same lines that if the  construction is done with different speeds for the different approximations of $X^k$ by
$X^k_n$, 
for example replacing $X^k_n$ by $\widetilde{X}^k_n$ defined by 
$$\widetilde{X}^k_n=X^k_{\ell_k(n)}\mbox{ with }\frac{\alpha_n}{\alpha_{\ell_k(n)}}\rightarrow a_k,$$
we would have
$$\alpha_n\mathbb{E}[(\varphi(\widetilde{X}^k_n)-\varphi(X^k))^2]\rightarrow a_k\int_{\mathbb{R}}\varphi^{\prime 2}\,d{\cal N}(0,1)$$
and we would obtain the generalized Mehler type error structure on the abstract Wiener space defined by $J$ associated with the quadratic form
$$e[f]=\sum_q a_q f_q^2$$
i.e. associated with the semi-group
$$p_t f=\sum_q<f,\xi_q>e^{-a_qt}\xi_q$$
(cf. [5] p113 et seq).\\

\noindent{\bf II.9. Approximation of a Poisson point process.}\\

Let $X$ be a r.v. with values in a metric space $E$ endowed with its Borel $\sigma$-field ${\cal F}$.

Let $X_n$ be an approximation of $X$ satisfying hypothesis (H3) with the sequence $\alpha_n$ and an algebra ${\cal D}_0$ composed of bounded
continuous functions (containing the constants and dense in $L^2(\mathbb{P}_X)$). We suppose that 
the Dirichlet form defined by $$\lim_n\alpha_n\mathbb{E}[(\varphi(X_n)-\varphi(X))^2]\quad\varphi\in{\cal D}_0$$ is {\it local}.
We denote by $\gamma[.]$ its square field operator.

Let $\mu$ be the law of $X$ on $(E,{\cal F})$. Let  $(X^j_n,X^j)$ be i.i.d. 
copies of $(X_n,X)$ and let $J$ be an integer valued r.v. with Poisson law of parameter 1 independent of the sequence $(X^j_n,X^j)$.

We consider the Poisson point processes
$$N_n=\sum_{j=1}^J \delta_{X^j_n}\qquad N=\sum_{j=1}^J \delta_{X^j}$$
($\sum_1^0$ meaning zero). $N_n$ and $N$ are r.v. with values in the space of point measures ${\cal M}_p$ on $(E,{\cal F})$ equipped with 
the smallest $\sigma$-field making all maps $A\mapsto \nu(A)\;A\in{\cal F}$ measurable for $\nu\in{\cal M}_p$.

We consider the algebra
$${\cal D}={\cal L}\{\nu\in{\cal M}_p\mapsto e^{i\int \varphi d\nu}\quad\varphi\in{\cal D}_0\}$$
{\bf Lemma 3.} {\it ${\cal D}$ is dense in $L^2(\mathbb{P}_N)$.}\\

\noindent{\bf Proof.} By the chaos decomposition, it is enough to prove that the constants and the elements of $L^2(\mathbb{P}_N)$ of the form
$\int f_1\,d\widetilde{N}\cdots\int f_k\,d\widetilde{N}$ where the functions $f_1,\ldots,f_k$ are measurable bounded on 
$(E,{\cal F})$ and where $\widetilde{N}$ denotes $N-\mu$, may be approached by elements of ${\cal D}$. Since the constants are in ${\cal D}$, it suffices to
reach $\int f_1\,dN\cdots\int f_k\,dN$. Now $\lim_{\lambda\rightarrow 0}(e^{\lambda\int \varphi dN}-1)/\lambda=\int \varphi dN$ and this gives easily 
the lemma.\hfill$\diamond$\\

We study the approximation of $N$ by $N_n$ by looking at 
$E_n=\alpha_n\mathbb{E}[(\Phi(N_n)-\Phi(N))^2]$ for $\Phi(\nu)=\exp[i\int \varphi d\nu]$. We may write
$$\begin{array}{rl}
E_n&{\displaystyle=\alpha_n\mathbb{E}[\exp(2i\sum_{j=1}^J\varphi(X^j_n))-2\exp(i\sum_{j=1}^J(\varphi(X^j_n)+\varphi(X^j))+\exp(2i\sum_{j=1}^J\varphi(X^j))]}\\
&{\displaystyle=\alpha_n\sum_{p=0}^\infty \frac{e^{-1}}{p!}[(\mathbb{E} e^{2i\varphi(X^1_n)})^p-2(\mathbb{E} e^{i(\varphi(X^1_n)+\varphi(X^1))})^p
+(\mathbb{E} e^{2i\varphi(X^1)})^p]}\\
&{\displaystyle=\alpha_n\sum_{p=0}^\infty \frac{e^{-1}}{p!}[(\mathbb{E} e^{2i\varphi(X^1_n)})^p-(\mathbb{E} e^{i(\varphi(X^1_n)+\varphi(X^1))})^p
+(\mathbb{E} e^{2i\varphi(X^1)})^p-(\mathbb{E} e^{i(\varphi(X^1_n)+\varphi(X^1))})^p]}\\
&{\displaystyle=\alpha_n\sum_{p=0}^\infty \frac{e^{-1}}{p!}[A_p+B_p]}
\end{array}
$$
We apply a similar idea to what we have done in section II.7 writing $A_p$ under the form
$$A_p=\mathbb{E}\left[e^{i\varphi(X^1_n)}(e^{i\varphi(X^1_n)}-e^{i\varphi(X^1)})\right]\sum_{k=1}^p
(\mathbb{E} e^{2i\varphi(X^1_n)})^{p-k}(\mathbb{E} e^{i(\varphi(X^1_n)+\varphi(X^1))})^{k-1}$$
Making now the additional assumption that $(X^1_n,X^1)\stackrel{d}{\Rightarrow}(X^1,X^1)$ and using the
 fact that $\varphi\in{\cal D}_0$ is continuous and bounded, we see that $\alpha_n(A_p+B_p)$ has the same limit as
$$\alpha_n\mathbb{E}[(e^{i\varphi(X^1_n)}-e^{i\varphi(X^1)})^2]p(\mathbb{E} e^{2i\varphi(X^1)})^{p-1}$$
but, since the form $\lim_n\alpha_n\mathbb{E}[(\varphi(X_n)-\varphi(X))^2]$ is local, theorem 2 on asymptotic error calculus applies and gives 
$$\lim_n\alpha_n\mathbb{E}[(e^{i\varphi(X^1_n)}-e^{i\varphi(X^1)})^2]=-\mathbb{E}[e^{2i\varphi(X^1)}\gamma[\varphi](X^1)]$$
and we obtain 
$$\begin{array}{rl}
\lim_n E_n&{\displaystyle=-\mathbb{E}[e^{2i\varphi(X)}\gamma[\varphi](X)]\sum_{p=0}^\infty\frac{e^{-1}}{p!}p(\mathbb{E} e^{2i\varphi(X)})^{p-1}}\\
&{\displaystyle=-e^{-1}e^{\mathbb{E}[e^{2i\varphi(X)}]}\mathbb{E}[e^{2i\varphi(X)}\gamma[\varphi](X)]}\\
&{\displaystyle=-\exp[-\int(1-e^{2i\varphi})d\mu]\int e^{2i\varphi}\gamma[\varphi]d\mu}.
\end{array}
$$
In order to recognize the limit, let us consider what we have called the white structure on $N$ associated with 
the structure $(E,{\cal F},\mu,{\cal D}_0,\gamma)$ (cf. [5] chapter VI \S3, cf. also [7]), its square field operator $\Gamma$ and its Dirichlet form
${\cal E}$ satisfy
$$
\begin{array}{l}
\Gamma[e^{i\int \varphi\,dN}]=-e^{2i\int\varphi\,dN}\int\gamma[\varphi]\,dN\\
{\cal E}[e^{i\int \varphi\,dN}]=-\frac{1}{2}\mathbb{E}[e^{2i\int\varphi\,dN}\int\gamma[\varphi]\,dN]
\end{array}
$$
which using the Laplace characteristic functional $\mathbb{E} e^{i\int fdN}=e^{-\int(1-e^{if})d\mu}$, may easily be seen to be equal to
$$-\frac{1}{2}\exp[-\int(1-e^{2i\varphi})d\mu]\int\gamma[\varphi]e^{2i\varphi}d\mu$$ what we obtained up to the factor 1/2. In conclusion\\

\noindent{\bf Proposition 15. }{\it Let $X_n$ be an approximation of $X$ satisfying } (H3) {\it on an algebra of continuous functions with a local 
asymptotic Dirichlet form and square field operator $\gamma$. Assuming in addition the weak convergence $(X_n,X)\stackrel{d}{\Rightarrow}(X,X)$
then the approximation $N_n$ of the Poisson point process $N$ constructed above satisfies} (H3) {\it with the same $\alpha_n$ and with asymptotic
Dirichlet form the so-called white form characterized by its square field operator}
$$\Gamma[\int fdN]=\int\gamma[f]dN.$$

\noindent{\bf II.9. Stochastic integral.}\\

We now consider a stochastic integral
$$Y=\int_0^1 H_s\,dB_s$$approximated by the sum
$$Y_n=\sum_{k=0}^{n-1}H_{\frac{k}{n}}(B_{\frac{k+1}{n}}-B_{\frac{k}{n}})$$
$(B_t)$ is a standard Brownian motion defined as the coordinate process of ${\cal C}([0,1])$ equipped with the Wiener measure, and $H_s=H_0+\int_0^s\xi_u\,dB_u+\int_0^s\eta_u\,du$
is an Ito process defined on the same space, processes $\xi$ and $\eta$ are adapted and regular in Malliavin sense. We suppose they satisfy 
$\sup_t\mathbb{E}[|\xi_t|^p+|\eta_i|^p]<+\infty$ for some $p>2$ and we will state their other regularity properties along the calculation.

In order to obtain the limit expressions we are looking for, we will use several times the integration by part formula
$$\mathbb{E}[u\delta U]=\mathbb{E}[<Du,U>_{\cal H}]$$
(cf.  for the notation [5] formula (15) p81). This technique has been already used with success by 
Clement, Kohatsu-Higa and Lamberton [11] to compute, for s.d.e. possibly with delay, an estimate of 
$\mathbb{E}[\varphi(Y_n)-\varphi(Y)]$, i.e. with our notation, an estimate of $<\overline{A}[\varphi],1>$. Let us note that this expression which is always equal to
$<\Abar[\varphi],1>$ erases the diffusive part of the bias and, since $\Abar$ is here a first order operator, as we will see in a moment,
this expression writes
$<\overline{A}[\varphi],1>=<\Abar[\varphi],1]=\mathbb{E}_Y[F\varphi^\prime]$ and reduces, when regularity allows an integration by parts, to the form
$\mathbb{E}[G\varphi]$. In the case of Ito type s.d.e. under rather general hypotheses, $\mathbb{E}[\varphi(Y_n)-\varphi(Y)]$ keeps the same order of magnitude
for the speed of convergence even when $\varphi$ is only bounded and measurable [2].

In this section we attempt to explicit the four bias operators for the above approximation problem. 
They occur with the sequence $\alpha_n=n$.\\

\noindent{\bf a) The local property is satisfied.}\\

\noindent{\bf Lemma 4. } {\it If $\xi$ and $\eta$ satisfy $\sup_t\mathbb{E}[|\xi_t|^p+|\eta_i|^p]<+\infty$ for some $p>2$ then 
$$n\mathbb{E}[|Y_n-Y|^{2+\alpha}]\rightarrow 0\quad\forall\alpha \;:\;2<2+\alpha\leq p$$}

\noindent{\bf Proof.} Let $\alpha$ be s.t. $2<2+\alpha\leq p$, by Burkholder-Davis-Gundy inequality
$$\mathbb{E}[|Y_n-Y|^{2+\alpha}]\leq C_1\mathbb{E}\int_0^1|H_s-H_{\frac{[ns]}{n}}|^{2+\alpha}\,ds.$$
Now 
$$\begin{array}{rl}
\|H_s-H_{\frac{[ns]}{n}}\|_{2+\alpha}&\leq (C_1\mathbb{E}[\int_{\frac{[ns]}{n}}^s\xi^2_udu)^{\frac{p}{2}}])^{\frac{1}{p}}
+(\mathbb{E}[(\int_{\frac{[ns]}{n}}^s|\eta_u|du)^p])^{\frac{1}{p}}\\
&\leq (s-\frac{[ns]}{n})^{\frac{1}{2}}(C_1\sup_t\mathbb{E}|\xi_t|^p)^{\frac{1}{p}}+(s-\frac{[ns]}{n})(\sup_t\mathbb{E}|\eta_t|^p)^{\frac{1}{p}}\\
\mathbb{E}[|H_s-H_{\frac{[ns]}{n}}|^{2+\alpha}]&\leq (s-\frac{[ns]}{n})^{\frac{2+\alpha}{2}}(C_2+o(1)).
\end{array}
$$
Hence $n\mathbb{E}[|Y_n-Y|^{2+\alpha}]\leq \frac{1}{n^{\alpha/2}}(C_3+o(1))$. Q.E.D.\hfill$\diamond$\\

It follows that if our test functions algebra ${\cal D}$ consists of bounded ${\cal C}^2$-functions with bounded derivatives, we have
$\lim_n n\mathbb{E}[|\varphi(Y_n)-\varphi(Y)|^{2+\alpha}]=0$ so that by remark 3 if we succeed in proving assumption (H3) with $\alpha_n=n$, the 
asymptotic Dirichlet form will be local.\\

\noindent{\bf b)} It follows also if we assume a little bit more for instance that the functions in ${\cal D}$ are ${\cal C}^3$ bounded with 
bounded derivatives, that in the study of $\overline{A}$, the expression $$n\mathbb{E}[(\varphi(Y_n)-\varphi(Y))\chi(Y)]$$ has the same limit as
$$n\mathbb{E}[(Y_n-Y)\varphi^\prime(Y)\chi(Y)+\frac{1}{2}(Y_n-Y)^2\varphi^{\prime\prime}(Y)\chi(Y)].$$
Similarly, in the study of $\widetilde{A}$, the expression 
$n\mathbb{E}[(\varphi(Y_n)-\varphi(Y))^2]$ has the same limit as $n\mathbb{E}[(Y_n-Y)^2\varphi^{\prime 2}(Y)]$.

For simplicity we shall suppose that the functions in ${\cal D}$ are ${\cal C}^\infty$ bounded with bounded derivatives.\\

\noindent{\bf c) Study of the symmetric bias operator.}

Let us remark first that in the study of $n\mathbb{E}[(Y_n-Y)^2\varphi^{\prime 2}(Y)]$ we may suppose $\eta\equiv 0$ 
and that $H$ be of the form $H_t=H_0+\int_0^t\xi_sdB_s$.

Indeed, putting $K_s=\int_0^s\eta_udu$ we have
$$ n\mathbb{E}[(\int_0^1(K_s-K_{\frac{[ns]}{n}})dB_s)^2\varphi^{\prime 2}(Y)]\leq Cn\int_0^1\mathbb{E}[(\int_\frac{[ns]}{n}^s\eta_udu)^2]ds=O(\frac{1}{n}).$$

Now for studying $n\mathbb{E}[(Y_n-Y)^2\varphi^{\prime 2}(Y)]$ we apply Ito's formula to the continuous martingale $Y_n-Y$ :
$$\begin{array}{rl}
(Y_n-Y)^2&=\int_0^12\int_0^t(H_s-H_{\frac{[ns]}{n}})dB_s(H_t-H_{\frac{[nt]}{n}})dB_t+\int_0^1(H_s-H_{\frac{[ns]}{n}})^2ds\\
&=(1)+(2)
\end{array}
$$

\noindent i) Contribution due to the second term.

Let us apply once more Ito's formula
$$\begin{array}{rl}
(\int_{\frac{[ns]}{n}}^s\xi_udB_u)^2&=2\int_{\frac{[ns]}{n}}^s
\int_{\frac{[nt]}{n}}^t\xi_udB_u\xi_tdB_t+\int_{\frac{[ns]}{n}}^s\xi^2_udu\\
&=(2,1)+(2,2).
\end{array}
$$
The contribution of the term (2,2) is $n\mathbb{E}[\int_0^1\int_{\frac{[ns]}{n}}^s\xi^2_ududs\varphi^{\prime 2}(Y)]$ which tends
 to $\frac{1}{2}\mathbb{E}[\int_0^1\xi^2_sds\varphi^{\prime 2}(Y)]$.

The contribution of the term (2,1) is zero. Indeed by integration by parts it is the limit of 
$$2n\int_0^1\mathbb{E}\int_{\frac{[ns]}{n}}^s\int_{\frac{[nt]}{n}}^t\xi_udB_u\xi_tD_t[\varphi^{\prime 2}(Y)]dtds$$
which by an other integration by part in order to get rid of the stochastic integral, gives 
$$2n\int_0^1\int_{\frac{[ns]}{n}}^s\int_{\frac{[nt]}{n}}^t\mathbb{E}[\xi_uD_u[\xi_tD_t[\varphi^{\prime 2}(Y)]]]dudtds$$
and is $O(\frac{1}{n})$as soon as $\mathbb{E}[\xi_uD_u[\xi_tD_t[\varphi^{\prime 2}(Y)]]]$ remains bounded.\\

\noindent ii) Contribution of the first term.

We shall show by several integration by parts that this contribution is zero. We are concerned by the limit of 
$$\begin{array}{l}
2n\mathbb{E}\int_0^1\int_0^t(H_s-H_{\frac{[ns]}{n}})dB_s(H_t-H_{\frac{[nt]}{n}})D_t[\varphi^{\prime 2}(Y)]dt\\
\quad\quad=2n\mathbb{E}\int_0^1\int_0^t(H_s-H_{\frac{[ns]}{n}})(H_t-H_{\frac{[nt]}{n}})D_sD_t[\varphi^{\prime 2}(Y)]dsdt\\
\quad\qquad+2n\mathbb{E}\int_0^1\int_0^t(H_s-H_{\frac{[ns]}{n}})D_s[H_t-H_{\frac{[nt]}{n}}]D_t[\varphi^{\prime 2}(Y)]dsdt\\
\qquad=(1,1)+(1,2).
\end{array}
$$
The term (1,1) decomposes in (1,11)+(1,12) with
$$\begin{array}{l}
(1,11) =2n\int_0^1\int_0^t\int_{\frac{[ns]}{n}}^s\mathbb{E}[\xi_u(H_t-H_{\frac{[nt]}{n}})D_uD_sD_t[\varphi^{\prime 2}(Y)]]dudsdt\\
(1,12) =2n\int_0^1\int_0^t\int_{\frac{[ns]}{n}}^s\mathbb{E}[\xi_uD_u[H_t-H_{\frac{[nt]}{n}}]D_sD_t[\varphi^{\prime 2}(Y)]]dudsdt
\end{array}
$$
we have 
$$
(1,11)=2n\int_0^1\int_0^t\int_{\frac{[ns]}{n}}^s\int_{\frac{[nt]}{n}}^t\mathbb{E}[\xi_vD_v[\xi_uD_uD_sD_t[\varphi^{\prime 2}(Y)]dvdudsdt$$
Thus (1,11) is $O(\frac{1}{n})$ as soon as the expectation inside is bounded. About (1,12) we get still two terms 
$$
(1,12)=2n\int_0^1\int_0^t\int_{\frac{[ns]}{n}}^s\mathbb{E}[\xi_u(\int_{\frac{[nt]}{n}}^tD_u[\xi_\alpha]dB_\alpha+
\xi_u1_{(\frac{[nt]}{n},t)}(u))D_sD_t
[\varphi^{\prime 2}(Y)]]dudsdt
$$
the second one requires $\frac{[ns]}{n}=\frac{[nt]}{n}$ and tends to zero. The first one may be written
$$2n\int_0^1\int_0^t\int_{\frac{[ns]}{n}}^s\int_{\frac{[nt]}{n}}^tD_u[\xi_\alpha]D_\alpha D_sD_t[\varphi^{\prime 2}(Y)]
d\alpha dudsdt=O(\frac{1}{n})$$

\noindent Let us come to the term (1,2) which may be written 
$$2n\int_0^1\int_0^t\mathbb{E}\left[(H_t-H_{\frac{[nt]}{n}})[\int_{\frac{[nt]}{n}}^tD_s[\xi_\alpha]dB_\alpha+\xi_s1_{(\frac{[nt]}{n},t)}(s)]
D_t[\varphi^{\prime 2}(Y)]\right]dsdt$$
still two terms, the second one requires $\frac{[ns]}{n}=\frac{[nt]}{n}$ and tends to zero. The first one is
$$
2n\int_0^1\int_0^t\mathbb{E}\left[\int_{\frac{[ns]}{n}}^s\xi_\beta dB_\beta\int_{\frac{[nt]}{n}}^t D_s[\xi_\alpha] dB_\alpha D_t[\varphi^{\prime 2}(Y)]\right]
dsdt$$
which may be handled as term (1,1).

Eventually, supposing $\xi$ and $\eta$ are bounded with bounded Malliavin derivatives up to order four, then 
$$n\mathbb{E}[(\varphi(Y_n)-\varphi(Y))^2]\rightarrow \frac{1}{2}\mathbb{E}[\int_0^1\xi^2_s\,ds\varphi^{\prime 2}(Y)].$$

\noindent{\bf d) Study of the theoretical bias $\overline{A}$.}

By the remark done in b) above, it is sufficient to study the limit of 
$$n\mathbb{E}[(Y_n-Y)\varphi^\prime(Y)\chi(Y)+\frac{1}{2}(Y_n-Y)^2\varphi^{\prime\prime}(Y)\chi(Y)].$$
The work is already done for the second term, it suffices to replace $\varphi^{\prime 2}(Y)$ by $\varphi^{\prime\prime}(Y)\chi(Y)$ in point c).

The first term may be written
$$
\begin{array}{l}
n\mathbb{E}[\int_0^1(H_{\frac{[ns]}{n}}-H_s)dB_s\varphi^\prime(Y)\chi(Y)]=n\mathbb{E}\int_0^1(H_{\frac{[ns]}{n}}-H_s)
D_s[\varphi^\prime(Y)\chi(Y)]ds\\
\qquad=n\mathbb{E}\int_0^1(\int_{\frac{[ns]}{n}}^s\xi_udB_u+\int_{\frac{[ns]}{n}}^s\eta_udu)D_s[\varphi^\prime(Y)\chi(Y)]ds\\
\qquad=(a)+(b)
\end{array}
$$
$(a)=n\int_0^1\int_{\frac{[ns]}{n}}^s\mathbb{E}[\xi_uD_uD_s[\varphi^\prime(Y)\chi(Y)]]duds$
which tends to $\frac{1}{2}\int_0^1\mathbb{E}[\xi_sD_sD_s[\varphi^\prime(Y)\chi(Y)]]ds$.

$$(b)=\int_0^1\int_{\frac{[ns]}{n}}^s\mathbb{E}[\eta_sD_s[\varphi^\prime(Y)\chi(Y)]duds\rightarrow 
\frac{1}{2}\int_0^1\mathbb{E}[\eta_sD_s[\varphi^\prime(Y)\chi(Y)]]ds$$
therefore we have 
\begin{equation}
\begin{array}{l}
\lim_n\mathbb{E}[(Y_n-Y)\varphi^\prime(Y)\chi(Y)]\qquad\qquad\\
\qquad\qquad=\frac{1}{2}\int_0^1\mathbb{E}[\xi_sD_sD_s[\varphi^\prime(Y)\chi(Y)]]ds
+\frac{1}{2}\int_0^1\mathbb{E}[\eta_sD_s[\varphi^\prime(Y)\chi(Y)]]ds.
\end{array}
\end{equation}
Hence, with the same hypotheses as for c), the conclusion is  :
$$\begin{array}{l}
\lim_n n\mathbb{E}[(\varphi(Y_n)-\varphi(Y))\chi(Y)]\qquad\\

{\displaystyle\qquad=\frac{1}{2}\int_0^1\mathbb{E}\left[\xi_sD_sD_s[\varphi^\prime(Y)\chi(Y)]+
\eta_sD_s[\varphi^\prime(Y)\chi(Y)]\right]ds
+\frac{1}{4}\int_0^1\mathbb{E}[\xi_s^2\varphi^{\prime\prime}(Y)\chi(Y)]ds}.
\end{array}
$$

\noindent{\bf e) Interpretation of the results.}

If we put $\frac{1}{2}\mathbb{E}[\int_0^1\xi_s^2 ds|Y=y]=\rho(y)$ the form ${\cal E}[\varphi]=\frac{1}{2}\int \rho\varphi^{\prime 2}\,d\mathbb{P}_Y$ is closable
iff the measure $\frac{1}{2}\rho\mathbb{P}_Y$ has a density satisfying the  Hamza condition (cf. [13] p105). In that case
the operator $\widetilde{A}$ exists and is uniquely defined by
$$<\widetilde{A}[\varphi],\chi>=-\frac{1}{4}\mathbb{E}[\int_0^1 \xi_s^2\,ds\varphi^\prime(Y)\chi^\prime(Y)].$$

\noindent The operator $\overline{A}$ would be defined by 
$$\begin{array}{rl}
<\overline{A}[\varphi],\chi>&=\lim_n n\mathbb{E}[(\varphi(Y_n)-\varphi(Y))\chi(Y)]\\
&=\frac{1}{2}\mathbb{E}[\int_0^1\xi_sD_sD_s[\varphi^\prime(Y)\chi(Y)]ds]+\frac{1}{2}\mathbb{E}[\int_0^1\eta_sD_s[\varphi^\prime(Y)\chi(Y)]ds]\\
&\quad +\frac{1}{4}\mathbb{E}[\int_0^1\xi_s^2(\varphi^\prime\chi)^\prime(Y) ds]-\frac{1}{4}\mathbb{E}[\int_0^1\xi_s^2\varphi^\prime(Y)\chi^\prime(Y) ds]
\end{array}
$$
provided that the righthand side may be put in the form of the lefthand side.

Sufficient conditions are easy to be listed.

(j) By the functional calculus the first term may be written
$$\frac{1}{2}\mathbb{E}\int_0^1\xi_s[(D_s[Y])^2(\varphi^\prime\chi)^{\prime\prime}(Y)+D_sD_s[Y](\varphi^\prime\chi)^\prime(Y)]ds$$
it will have the desired form as soon as the measures $\int_0^1\xi_s(D_s[Y])^2ds.\mathbb{P}_Y$ and \break
$\int_0^1\xi_sD_sD_s[Y]ds.\mathbb{P}_Y$ will be sufficiently regular to allow an integration by parts.

(jj) The second term is equal to
$$\frac{1}{2}\mathbb{E}[\int_0^1 \eta_sdB_s\varphi^\prime(Y)\chi(Y)].$$

(jjj) The third term requires the conditions of an integration by part.

(jv) The fourth term is $<\widetilde{A}[\varphi],\chi>$.\\

We see  that the operator $\widetilde{A}$ exists under a quite simple condition involving only the law of the pair $(\int_0^1\xi_s^2ds, Y)$.
The regularity conditions insuring the existence of the operators $\overline{A}$ or $\underline{A}$ are more intricate. When they hold, putting
\begin{equation}
<\Abar[\varphi],\chi>=\frac{1}{2}\mathbb{E}\int_0^1(\xi_sD_sD_s+\eta_sD_s)[\varphi^\prime(Y)\chi(Y)]ds+\frac{1}{4}\mathbb{E}[\int_0^1\xi^2_sds(\varphi^\prime\chi)^\prime(Y)]
\end{equation}
we have $\overline{A}=\Abar+\widetilde{A}$ and $\underline{A}=-\Abar+\widetilde{A}$ and by the 
general theory, hypotheses (H1) to (H3) being fulfilled and the Dirichlet form being local, $\Abar$ is a first order operator, 
as may be seen also on the obtained form (12).\\

\noindent{\bf Remark 11. } Our approach is direct. But the heaviest part of the proof i.e. the proof of 
$$n\mathbb{E}[(\varphi(Y_n)-\varphi(Y))^2]\rightarrow \frac{1}{2}\mathbb{E}\int_0^1\xi^2_sds\varphi^{\prime 2}(Y),$$
may be considerably  shortened if we use a result of weak convergence like
$$(\sqrt{n}(Y_n-Y),Y)\stackrel{d}{\Rightarrow}(\frac{1}{\sqrt{2}}\int_0^1\xi_sdW_s,Y)$$
with an ``extra" Brownian motion $W$ independent of $Y$ and $\xi$.

This gives our results thanks to the uniform integrability of $n(Y_n-Y)^2$ which is a consequence of the inequality 
$\mathbb{E}[|Y_n-Y|^{2+\alpha}]\leq n^{-\frac{2+\alpha}{2}}(C_3+o(1))$ established in the proof of lemma 4.

Such a weak convergence result has been obtained long time ago by Rootzen  [28 ]  for the case where the process $H$ has the form 
$H_s=f(B_s,s)$. This kind of weak convergence results for stochastic integrals have been now considerably extended, 
see especially
[16], [18], [29], [19], [21], [21],  [17],  [15]. Our approach to s.d.e. in the next section 
is based on such results.\\

\noindent{\bf II.10. Stochastic differential equations and Euler scheme.}\\

As we have just explained we will base our approach on  results on  convergence in law, in particular on the article of Jacod and Protter 
[17]. We consider only the case of a continuous semi-martingale in which the main ideas already appear.

Let $X=(X^i)_{i=1,\ldots,d}$ be a continuous semi-martingale with values in $\mathbb{R}^d$ vanishing at zero defined on the stochastic basis 
$(\Omega, {\cal F},({\cal F}_t),\mathbb{P})$. For $t\in[0,1]$ we consider the $q$-dimensional s.d.e.
\begin{equation}
dY_t=f(Y_t)\,dX_t\qquad Y_0=y_0
\end{equation}
where $y_0\in\mathbb{R}^q$, $f$ is ${\cal C}^1$ from $\mathbb{R}^q$ into $\mathbb{R}^{q\times d}$ with at most linear growth 
($|f(x)|\leq K(1+|x])$ denoting $|.|$ the norms on $\mathbb{R}^k$). It is known that (13) has a unique strong solution. We study the resolution of (13) 
by the Euler scheme :
$$
dY^n_t=f(Y^n_{\frac{[nt]}{n}})\,dX_t\qquad Y_0^n=y_0
$$
where $[nt]$ is the entire part of $nt$.

We denote $U^n_t=Y^n_t-Y_t$ the error process. $U^n$ as process with values in ${\cal C}([0,1])$ tends to zero in probability (as soon as $f$
is locally Lipschitz with at most linear growth [17]).

It is supposed that $X=M+A$ where $M$ is a continuous local martingale vanishing at zero with values in $\mathbb{R}^d$ and $A$ is a continuous
finite variation adapted process vanishing at zero satisfying
$$\begin{array}{l}
A^i_t=\int_0^ta^i_s\,ds\;\mbox{ with }\int_0^1(a^i_s)^2\,ds<+\infty\mbox{ a.s.}\\
<M^i,M^j>_t=\int_0^t c^{ij}_s\,ds\;\mbox{ with }\int_0^1(c^{ij}_s)^2ds<+\infty\mbox{ a.s.}
\end{array}
$$
then for every starting point $y_0$ and for all function $f$ ${\cal C}^1$ with at most linear growth, the process $\sqrt{n}U^n$ converges in law
on ${\cal C}([0,1])$ to the solution to 
$$
dU^i_t=\sum_{j=1}^d\sum_{k=1}^q\frac{\partial f^{ij}}{\partial x_k}(Y_t)\left[U^k_t\,dX^j_t-\sum_{\ell=1}^d f^{k\ell}(Y_t)\,dZ^{\ell j}_t\right],\qquad U^i_0=0,
$$
 $Z$ being given by 
$$
Z^{ij}_t=\frac{1}{\sqrt{2}}\sum_{k,\ell=1}^q\int_0^t\sigma^{ik}_s\sigma^{j\ell}_s\,dW^{k\ell}_s
$$
where $W$ is a standard $q^2$-dimensional Brownian motion defined on an extension of the space  independent of $X$ and 
$\sigma$ is a matrix of processes s.t. $(\sigma\sigma^t)^{ij}=c^{ij}$ which exists as soon as $q\geq d$ 
case to which the question may be always reduced.

The proof consists of the conjonction of theorems 3.3, 5.1 and 5.5 of [17] and their proofs.\\

In order to study the hypotheses (H1) to (H3) we consider the algebra ${\cal D}$ of the linear conbinations
of functions $\varphi$ defined on ${\cal C}([0,1])$ by
$$\varphi(Y)=e^{i<u_1,Y_{t_1}>+\cdots+i<u_r,Y_{t_r}>}\qquad u_\ell\in\mathbb{R}^q\qquad t_\ell\in[0,1]\qquad \ell=1,\ldots,r$$
and the sequence $\alpha_n=n$.\\

\noindent{\bf a) Symmetric bias operator.}

We study $n\mathbb{E}[(\varphi(Y^n)-\varphi(Y))^2]$.\\

\noindent{\bf Lemma 5. }{\it If for fixed $t$ the sequence $n|Y^n_t-Y_t|^2=|\sqrt{n}U^n_t|^2$ is uniformly integrable, }
\begin{equation}
n\mathbb{E}[(\varphi(Y^n)-\varphi(Y))^2]\rightarrow 
\mathbb{E}\left[\left(\sum_{j=1}^q\sum_{\ell=1}^r U^j_{t_\ell}\frac{\partial \varphi}{\partial y^j_{t_\ell}}(Y)\right)^2\right].
\end{equation}

\noindent{\bf Proof.} Let us argue in the case $q=r=1$, the general case being similar.
$$\begin{array}{l}
n\mathbb{E}[(\varphi(Y^n)-\varphi(Y))^2]=\mathbb{E}[n(Y^n_t-Y_t)^2(\int_0^1\varphi^\prime(Y_t+\lambda(Y^n_t-Y_t))d\lambda)^2]\qquad\qquad\\
\qquad\qquad\leq \mathbb{E}[\{n(Y^n_t-Y_t)^2-(n(Y^n_t-Y_t)^2)\wedge a\}\|\varphi^\prime\|^2_\infty]\\
\qquad\qquad\qquad+\mathbb{E}[\{(n(Y^n_t-Y_t)^2)\wedge a\}(\int_0^1\varphi^\prime(Y_t+\lambda(Y^n_t-Y_t))d\lambda)^2]
\end {array}
$$
By the uniform integrability the first term may be made smaller than $\varepsilon>0$ uniformly in $n$ by suitable choice of $a$, then the second term
goes to zero by the weak convergence of $n(Y^n-Y)$ and the convergence in probability of $Y^n_t-Y_y$ to zero.\hfill$\diamond$\\

\noindent{\bf Remark 12. } In the classical case of an s.d.e. defining a diffusion process from a Brownian motion, if the coefficients are regular, 
for instance ${\cal C}^\infty$ with bounded derivatives, it is known that $\sqrt{n}\|Y^n_t-Y_t\|_p$ is bounded for any $p\in[1,+\infty[$,
the uniform integrability of $n|Y^n_t-Y_t|^2$ follows.\\

Considering that $X$ and $W$ are defined on a product space whose samples are denoted $\omega$ and $\hat{\omega}$, formula (14)
shows that if hypothesis (H3) is verified and if $n|U^n_t|^2$ is uniformly integrable, the limit Dirichlet form satisfies $Y_t\in\mathbb{D}$ and its square field
operator satisfies
$$\Gamma[Y^j_t]=\hat{\mathbb{E}}[(U^j_t)^2].$$
In other words, the limit process $U(\omega,\hat{\omega})$ appears to be a gradient in the sense of Dirichlet forms of 
the process $Y$ : we may write
\begin{equation}
(Y_t)^\#(\omega,\hat{\omega})=U_t(\omega,\hat{\omega})
\end{equation}
and formula (14) follows by the chain rule.\\

The remaining question is whether the form defined on ${\cal D}$ by (14) is closable in $L^2({\cal C}([0,1]),\mathbb{P}_Y)$. To this 
question we have yet only an answer in the simplest case where $q=1$. When
$$dY_t=a(Y_t,t)dB_t+b(Y_t,t)dt$$
with $a,b$ ${\cal C}^1$ with at most linear growth, the process $U$ is given by 
$$U_t=N_t\int_0^t\frac{a(Y_s,s)a^\prime_y(Y_s,s)}{\sqrt{2}N_s}\,dW_s$$
with
$$N_t=\exp\{\int_0^t a^\prime_y(Y_s,s)dB_s-\frac{1}{2}\int_0^ta^{\prime 2}_y(Y_s,s)ds+\int_0^t b^\prime_y(Y_s,s)ds\}.$$
Let us denote $({\cal E}^\theta_{ou},\mathbb{D}^\theta_{ou})$ the Dirichlet form on the Wiener space of type Ornstein-Uhlenbeck with deterministic weight $\theta$, 
and let us denote $D^\theta_{ou}$ its gradient operator defined with the auxiliary Hilbert space $L^2([0,1],dt)$. We have\\

\noindent{\bf Proposition 16. }{\it If the coefficient $a$ satisfies $\mathbb{E}\int_0^1a^{\prime 2}_y(Y_s,s)ds<+\infty$ and if 
$a^{\prime 2}_y(Y_s,s)\geq\theta(s)>0$, hypothesis} (H3) {\it is fulfilled. The asymptotic Dirichlet form is the image by $Y$ of the
 form $({\cal E}_w,\mathbb{D}_w)$ defined on the Wiener space by }
$$\mathbb{D}_w=\{F\in\mathbb{D}^\theta_{ou}\;:\;\int_0^1\mathbb{E}[(D^\theta_{ou}[F](t))^2\frac{a^{\prime 2}_y(Y_t,t)}{\theta(t)}]dt<+\infty\}$$
$${\cal E}_w[F]=\frac{1}{4}\int_0^1\mathbb{E}[(D^\theta_{ou}[F](t))^2\frac{a^{\prime 2}_y(Y_t,t)}{\theta(t)}]dt.$$

\noindent The proof has been exposed at the Fifth Seminar on Stochastic Analysis, Random Fields and Application at 
Ascona in 2005 and will appear in the proceedings.

The form $({\cal E}_w,\mathbb{D}_w)$ admits the square field operator 
$$\Gamma_w[F]=\frac{1}{2}\int_0^1(D^\theta_{ou}[F](t))^2\frac{a^{\prime 2}_y(Y_t,t)}{\theta(t)}dt.$$
Putting $\xi_t=\frac{1}{2}a^{\prime 2}_y(Y_t,t)$ the operator $\widetilde{A}$ is given by
$$\widetilde{A}[\varphi](y)=\mathbb{E}[A_w[\varphi(Y)]|Y=y]$$
where $A_w[\varphi(Y)]=-\frac{1}{2}\delta^\theta_{ou}[\frac{\xi}{\theta}D^\theta_{ou}[F]]$, and $\delta^\theta_{ou}$ being the Skorokod stochastic integral operator
associated with $({\cal E}^\theta_{ou},\mathbb{D}^\theta_{ou})$.

From the concrete point of view of error calculus, the relation 
\begin{equation}
Y_t^\#=U_t
\end{equation}
is the most important. It allows to propagate errors by the chain rule and  using also, in order to manage limit objects, the fact that 
the operator $\#$ is closed what is a consequence of the closedness of the form.\\

\noindent{\bf b) The theoretical bias operator.}

As in the case of the approximation of a stochastic integral (cf. section II.9) the operator $\overline{A}$ involves an iterated gradient.

The main part of the calculation has been performed by Malliavin and Thalmaier ([24] and [25]) and we adopt their hypotheses :
$Y$ is solution of the s.d.e.
$$dY_t=a(Y_t)dB_t+b(Y_t)dt$$
where $B$ is a $(d-1)$-dimensional Brownian motion and where the matrix $a$ and the function $b$ are ${\cal C}^\infty$ 
with bounded derivatives.

The operator $\overline{A}$ is given by $\lim_n n\mathbb{E}[(\varphi(Y^n)-\varphi(Y))\chi(Y)]$. Since ${\cal D}$
 consists of functions of finite dimensional marginals, we restrict for simplicity to marginals of order one 
and to the case where $B$ and $Y$ are scalar, we have (cf. lemma 4) :
$$\lim_n n\mathbb{E}[(\varphi(Y^n_t)-\varphi(Y_t))\chi(Y_t)]=\lim_n n\mathbb{E}[(Y^n_t-Y_t)\varphi^\prime(Y_t)\chi(Y_t)
+\frac{1}{2}(Y^n_t-Y_t)^2\varphi^{\prime\prime}(Y_t)\chi(Y_t)]$$
P. Malliavin and A. Thalmaier have computed the first term which may be pulled back on the Wiener space
\begin{equation}
\lim_n n\mathbb{E}[(Y^n_t-Y_t)\varphi^\prime(Y_t)\chi(Y_t)]=\int_0^1\mathbb{E}[a_{11}(Y_s)D_sD_sF+b_1(Y_s)D_sF+c_1(Y_s)F]ds
\end{equation}
where $F=\varphi^\prime(Y_t)\chi(Y_t)$ and where $a_{11},b_1,c_1$ are functions of the coefficients $a,b$ and 
of their four first derivatives. It should be noted the similarity between the above formula (17) and formula (11) obtained
 for a stochastic integral ((17) reduces to (11) when $b=0$).

The second term is consequence of the preceding results on convergence in law.
\begin{equation}
\frac{n}{2}\mathbb{E}[(Y^n_t-Y_t)^2\varphi^{\prime\prime}(Y_t)\chi(Y_t)]\rightarrow
\frac{1}{2}\mathbb{E}[U^2_t\varphi^{\prime\prime}(Y_t)\chi(Y_t)].
\end{equation}
We see, by formulae (17) and (18) that the operator $\overline{A}$ is the image by $Y$ of a singular 
distribution operator on the Wiener space.\\

We have to conclude that, up to now, the study of the approximation of the solution of an s.d.e. by the Euler 
scheme is far from being achieved : the operator $\overline{A}$ is yielded by the quoted recent works 
but the existence of the operator $\widetilde{A}$ (hence of the Dirichlet form) is only shown in a very particular case.\\

\noindent{\Large\sf III. Conclusive comments.}\\

We focuse in this conclusion on remarks concerning the comparison between deterministic and stochastic approximation.\\

\noindent{\it The hypothesis of uniqueness of the approximation of order $n$.}

Let us consider a situation where given $Y$ the approximation $Y_n$ is completely determined, i.e. for $\mathbb{P}_Y$-a.e. $y$, the conditional law of 
$Y_n$ given $Y=y$, has the form $\delta_{\eta_n(y)}$. We call this assumption of uniqueness hypothesis (U).\\

\noindent{\it Example.} Such a hypothesis if often {\it implicitely} supposed when numerical results are given under the form
\begin{equation}\begin{array}{l}
Y_3=2.3769\pm 10^{-4}\\
Y_5=2.376985\pm10^{-6}\\
Y_7=2.37698534\pm10^{-8}\\
\ldots
\end{array}
\end{equation}
and it is underlying the concept of {\it number of significant digits}.

Indeed, let us take the decimal representation of real numbers in $[0,1]$ :
$$y=\sum_{n=0}^\infty\frac{a_n}{10^{n+1}}\quad\mbox{with }\quad a_n\in\{0,1,\ldots, 9\}$$
If the $a_n$'s are drawn independently uniforly on $\{0,1,\ldots, 9\}$ the random variable $Y=\sum_{n=0}^\infty\frac{a_n}{10^{n+1}}$
is uniformly ditributed on $[0,1[$ and as soon as $y$ is not decimal, which is a negligeable set, the expansion of $y$ is unique, so that the above
hypothesis (U) is fulfilled for the approximation $Y_n=\sum_{k=0}^n\frac{a_k}{10^{k+1}}$.\\

Some martingales satisfy hypothesis (U). Let us consider $(\Omega,{\cal A},\mathbb{P})$ with an increasing sequence of sub-$\sigma$-fields ${\cal B}_n$ generated
by countable partitions ${\cal P}_n$ of $\Omega$. Then for $Y\in L^1$, $Y_n=\mathbb{E}[Y|{\cal B}_n]$ satisfies (U) with
$$\eta_n(y)=\sum_{A\in{\cal P}_n}1_A(y)\frac{\mathbb{E}[Y1_A]}{\mathbb{P}(A)}.$$
This happens in particular for Haar systems (cf. [26] chap. III \S3).\\

\noindent{\bf Proposition 17. } {\it Suppose hypothesis } (U). {\it If for  $\alpha_n\rightarrow+\infty$ and an algebra ${\cal D}$, }
$$\alpha_n\mathbb{E}[(\varphi(Y_n)-\varphi(Y))\chi(Y)]\rightarrow<\overline{A}[\varphi],\chi>_{\mathbb{P}_Y}\quad\quad
\forall \varphi\in{\cal D},\quad\forall \chi\in L^2(\mathbb{P}_Y),$$
 {\it then} (H1) {\it to } (H3) {\it hold, $\overline{A}=-\underline{A}=\Abar$ are first order operators and 
$\widetilde{A}=0$.}\\

\noindent{\bf Proof.}
The sequence $\alpha_n(\mathbb{E}[\varphi(Y_n)|Y=y]-\varphi(y))$ is weakly bounded in $L^2(\mathbb{P}_Y)$
 hence strongly bounded, i.e. 
$$\alpha_n^2\int(\mathbb{E}[\varphi(Y_n)|Y=y]-\varphi(y))^2\;\mathbb{P}_Y(dy)\leq K.$$
Now $\mathbb{E}[\varphi(Y_n)|Y=y]=\varphi(\eta_n(y))$, hence 
$$\alpha_n\mathbb{E}[(\varphi(Y_n)-\varphi(Y))^2]=\alpha_n\int(\varphi(\eta_n(y))-\varphi(y))^2\;\mathbb{P}_Y(dy)
\leq \frac{K}{\alpha_n}\rightarrow 0.$$
The Dirichlet form is zero, hence it is local and $\Abar$ is a first order operator.\hfill$\diamond$\\

Examples in part II show that in many probabilistic approximations, hypothesis (U) does
 not hold. The law of $Y_n$ given $Y=y$ has a non zero variance. Polya's urn is a generic example showing that information at each step 
cannot be resumed by boxes of size $\pm10^{-k}$ but by standard deviation of laws whose support doesn't go in general to zero.

When we are interested in computing a sample of a random quantity, we have to display the result with specifications adapted
 to the stochastic case.

The interest of such specifications is particularly clear in infinite dimension when we have
 to compute $\omega$ by $\omega$ a path of a process. For example in the GPS or GALILEO systems 
when modelling the ionosphere by a spatio-temporal process, computing a sample is necessary to obtain
 the shift in the signals coming from several satellites. The accuracy of this sample is important to get the accuracy
 of the whole positionning system.

For such numerical computations of sample paths, we suggest that, as much as possible, the following specifications be displayed :
(i) the law $\mathbb{P}_Y$ of $Y$, (ii) the sequence $\alpha_n$, (iii) the theoretical and practical bias operators 
$\overline{A}$ and $\underline{A}$.

Then the operator $\widetilde{A}$, the Dirichlet form and the square field operator follow and the 
approximation $Y_n(\omega)$ may be (if the Dirichlet form is local)  the starting point of an error calculus for the studied model.
Non locality of the form, when it happens, is also a precious warning to be particularly carefull in the sensitivity analysis.

\begin{list}{}
{\setlength{\itemsep}{0cm}\setlength{\leftmargin}{0.5cm}\setlength{\parsep}{0cm}\setlength{\listparindent}{-0.5cm}}
  \item\begin{center}
{\small REFERENCES}
\end{center}\vspace{0.4cm}

[1] {\sc Ancona, A.} ``Continuit\'e des contractions dans les espace de Dirichlet" in {\it S\'emi\-naire de Th\'eorie du Potentiel,
Paris $n^0$2}, p1-26, Lect. N. in Math 963, Springer 1976.

[2] {\sc Bally, V., Talay, D.} ``The Euler scheme for stochastic differential 
equations: error analysis with Malliavin calculus" {\it Math. and Computers in Simulation} 38, p35-41, (1995).

[3] {\sc Billingsley, P.} {\it Convergence of Probability Measures} Wiley (1968),

[4] {\sc Bouleau N.} ``Error calculus and path sensitivity in Financial models", 
{\it Mathematical Finance} vol 13/1, 115-134, (2003).

[5] {\sc Bouleau N.} {\it Error Calculus for Finance and Physics, the Language of Dirichlet Forms}, De Gruyter, 2003.

[6] {\sc Bouleau N.} ``Th\'eor\`eme de Donsker et formes de Dirichlet" {\it Bull. Sci. Math. }129, (2005), 369-380.

[7] {\sc Bouleau N.} ``Improving Monte Carlo simulations by Dirichlet forms" {\it C. R. Acad. Sci. Paris} Ser I (2005)

[8] {\sc Bouleau, N. } and {\sc Chorro, Chr.} ``Error structures and parameter estimation" {\it C. R. Acad. Sci. Paris }s\'er I 338, (2004) 305-310.

[9] {\sc Bouleau N., Hirsch F.}  {\it Dirichlet Forms and Analysis on Wiener Space,} De Gruyter, (1991).

[10] {\sc Chorro, Chr.} {\it Calculs d'erreur par formes de Dirichlet, liens avec l'information de Fisher et les th\'eor\`emes limites   } Th\`ese, Universit\'e Paris 1, (2005).

[11] {\sc Clement, I.; Kohatsu-Higa, A.; Lamberton, D.} ``A duality approach for the weak approximation of a stochastic differential equation"
{\it Annals of Applied Probability} (to appear).

[12] {\sc Dellacherie, Cl., Meyer, P.-A.,} {\it Probabilit\'es et Potentiel,} Hermann 1987.

[13] {\sc Fukushima, M.; Oshima, Y.; Takeda, M.} {\it Dirichlet forms and symmetric Markov processes}, De Gruyter 1994.

[14] {\sc Hall, P., Heyde, C. C., } {\it Martingale limit theory and its applications} Academic Press (1980).

[15] {\sc Hayashi, T., Mykland, P.A.} ``Evaluating hedging errors : an asymptotic approach" {\it Math. Finance}, vol 15, No 2, 309-343, (2005)

[16] {\sc Jacod, J.,} ``Th\'eor\`emes limites pour les processus" {\it Lect. Notes Math. } vol 1117, Springer 1985.

[17] {\sc Jacod, J., Protter, Ph.} ``Asymptotic error distributions for the Euler method for stochastic 
differential equations'' {\it Ann. Probab.} 26, 267-307, (1998)

[18] {\sc Jacod, J., Shiryaev, A.N.,} {\it Limit Theorems for Stochastic Processes,} Springer 1987.

[19] {\sc Jakubowski, A., M\'emin, J., Pag\`es, G.} ``Convergence en loi de suite d'int\'egrales stochastiques sur l'espace de Skorokhod" {\it 
Probab. Th. Rel. Fields} 81, 111-137, 1989.

[20] {\sc Kolmogorov, A.N.},  ``Ueber die analytischen Methoden in der Wahrscheinlichkeitsrechnung" {\it Math. Ann.}
104-108, (1931)

[21] {\sc Kurtz, Th; Protter, Ph.} ``Wong-Zakai corrections, random evolutions and simulation schemes for SDEs" {\it Stochastic Analysis} 331-346,
Acad. Press, 1991.

[22] {\sc Kurtz, Th; Protter, Ph.} ``Weak limit theorems for stochastic integrals ans stochastic differential equations" {\it Ann. Probab.} 19, 1035-1070, 1991.

[23] {\sc Ma, Z.-M., R\"{o}ckner, M.} {\it Introduction to the Theory of (Non-symmetric) Dirichlet Forms} Springer 1992.

[24] {\sc Malliavin, P., Thalmaier, A.} ``Numerical error for SDE: asymptotic espansion and hyperdistributions" {\it Note C. R. A. S.} sI, vol 336, n¡10, p851, 2003.

[25] {\sc Malliavin, P., Thalmaier, A.} {\it Stochastic Calculus of Variations in Mathematical Finance}, Springer 2005.

[26] {\sc Neveu, J.} {\it Martingales \`a temps discret} Masson (1972).

[27] {\sc Nualart, N.} : {\it The Malliavin calculus and related topics}. Springer, 1995.

[28 ] {\sc Rootz\'en, H.} ``Limit distribution for the error in approximation of stochastic intergrals'' {\it Ann. Probab. } 8, 241-251, (1980).

[29] {\sc Slomi\'nski, L.} ``Stability of strong solutions of stochastic differential equations'' {Stochastic Process. Appl.} 31, 173-202, (1989).

\end{list}

\end{document}